\documentclass[11pt, oneside]{article} 

\usepackage[utf8]{inputenc}
\usepackage[T1]{fontenc}
\usepackage{xcolor} 
\usepackage{color} 

\usepackage{amsmath}
\usepackage{amssymb} 
\usepackage{mathrsfs} 
\usepackage{bm} 
\usepackage{bbm} 

\usepackage{amsthm} 
\theoremstyle{plain}
\newtheorem{theorem}{Theorem} 

\newtheorem{example}[theorem]{Example}

\usepackage{algorithm2e}
\SetKwInput{KwInputs}{Inputs}
\usepackage[noend]{algpseudocode}
\makeatletter
\def\BState{\State\hskip-\ALG@thistlm}
\makeatother

\RequirePackage[numbers]{natbib}
\RequirePackage[colorlinks,citecolor=blue,urlcolor=blue]{hyperref}
\usepackage{graphicx} 
\graphicspath{{Figures/} }
\usepackage{subcaption}

\usepackage{tikz} 
\usetikzlibrary{positioning ,shadows,matrix}

\usepackage[english]{babel} 

\usepackage{hyperref} 
\hypersetup{
    pdfpagemode={UseOutlines},
    bookmarksopen,
    pdfstartview={FitH},
    colorlinks,
    linkcolor={blue}, 
    citecolor={blue}, 
    urlcolor={blue} 
}

\usepackage[letterpaper,margin=1.25in]{geometry} 

\usepackage{enumitem} 

\def\iid{\overset{\textnormal{iid}}{\sim}} 

\makeatletter  
\@ifundefined{@currsizeindex} 
  {\RequirePackage{relsize}\let\dolarger\relsize} 
  {\def\dolarger#1{\larger[#1]}} 
\newcommand*\@@bigtimes[2]{\vphantom{\prod} 
  \vcenter{\hbox{\dolarger{4}$\m@th#1\mkern-2mu\times\mkern-2mu$}}} 
\newcommand*\bigtimes{\mathop{\mathpalette\@@bigtimes\relax}\displaylimits} 
\makeatother

\def\E{\mathbb{E}}\def\H{\mathbb{H}}\def\I{\mathbb{I}}\def\N{\mathbb{N}}\def\R{\mathbb{R}}\def\1{\mathbbm{1}}\def\Cov{\mathbbm{C}\textnormal{ov}}

\def\Lcal{\mathcal{L}}\def\Ocal{\mathcal{O}}\def\Scal{\mathcal{S}}

\def\Enorm{\textnormal{E}}

\def\eps{\varepsilon}

\title{\bf A Bayesian approach with Gaussian priors to the inverse problem of source identification in elliptic PDEs}
\author{Matteo Giordano \\ \\ ESOMAS Department, University of Turin, \\
Corso Unione Sovietica 218 bis, Turin, Italy}
\date{} 
\begin{document}

\maketitle

\abstract{
We consider the statistical linear inverse problem of making inference on an  unknown source function in an elliptic partial differential equation from noisy observations of its solution. We employ nonparametric Bayesian procedures based on Gaussian priors, leading to convenient conjugate formulae for posterior inference. We review recent results providing theoretical guarantees on the quality of the resulting posterior-based estimation and uncertainty quantification, and we discuss the application of the theory to the important classes of Gaussian series priors defined on the Dirichlet-Laplacian eigenbasis and Matérn process priors. We provide an implementation of posterior inference for both classes of priors, and investigate its performance in a numerical simulation study. The reproducible code is available at: \url{https://github.com/MattGiord/Bayesian-Source-Identification
}

\medskip

}

\tableofcontents

\section{Introduction}
\label{Sec:Intro}

Linear inverse problems consist in the task of recovering unknown objects or physical quantities from linear indirect noisy measurements. A widespread mathematical formulation for these problems postulates that the recovery target be an element $f$ of a Hilbert space $\H_1$, and that the data arise according to the equation
\begin{equation}
\label{Eq:GenIP}
	Y^\eps=G(f)+\eps W,
\end{equation}
where $G:\H_1\to\H_2$ is a linear operator between $\H_1$ and another Hilbert space $\H_2$, $W$ is additive observational noise and $\eps>0$ is the noise level. In view of the central limit theorem, normality of the measurement errors can often be maintained, whereby $W$ is assumed to be a white noise process indexed by $\H_2$. The goal is then to estimate $f$ from an observed realisation of $Y^\eps$.

	Observation schemes as in equation \eqref{Eq:GenIP} are found in a variety of scientific fields and engineering applications, including medical imaging \cite{BP06}, geophysics \cite{ST99}, acoustics \cite{collins1994inverse} and finance \cite{baumeister2013inverse}. For example, Computerised Tomography (CT), a technique to obtain detailed images of the human body and information about the density variation of the tissues, is based on a mathematical model (related to the `Radon transform') for the absorption of $X$-rays. Similar concepts underpin many other medical imagining techniques, such as Magnetic Resonance (MR) and Positron Emission Tomography (PET); see  \cite{BP06} for further details and references.

	 In many such applications, the unknown $f$ may be characterised as a  functional coefficient governing a partial differential equation (PDE), while the observed object $G(f)$ is the corresponding PDE solution. The `forward operator' is then the coefficient-to-solution map $G:f\mapsto G(f)$. See the monograph \cite{I17} for an extensive overview on PDE-based inverse problems.

	In the present article, we shall mostly focus on the representative example of `source identification' in elliptic PDEs. Many of the ideas developed hereafter have a natural application to other linear inverse problems. Let $\Ocal\subset \R^d$ be an open and bounded set with smooth boundary $\partial\Ocal$. Let the unknown (square-integrable) function $f\in \H_1\equiv L^2(\Ocal)$ be the `source term' in the elliptic PDE with zero Dirichlet boundary conditions,
\begin{equation}
\label{Eq:PDE}
\begin{split}
	\nabla\cdot(c\nabla u)-f&=0,\ \ \text{on}\ \ \Ocal,\\
	u&=0, \ \ \text{on}\ \ \partial\Ocal,
\end{split}
\end{equation}
where $\nabla\cdot$ and $\nabla$ denote, respectively, the divergence and gradient operators, and where the smooth and positive `diffusion coefficient' $c\in C^\infty(\overline\Ocal)$, $\inf_{x\in\Ocal}c(x)>0$, is assumed to be known. By standard elliptic theory, e.g.~\cite[Chapter 2]{LM72} and \cite[Chapter 6]{E10}, for any $f\in L^2(\Ocal)$ there exists a unique (weak) solution $G(f)\equiv u$ in the Sobolev space $H^1(\Ocal)\subset L^2(\Ocal)\equiv \H_2$, giving rise to a linear (injective, self-adjoint and compact) operator $G:L^2(\Ocal)\to L^2(\Ocal)$; see Appendix B in \cite{GK20} for further details. We then consider the problem of estimating the source $f$ from observations $Y^\eps$ arising as in \eqref{Eq:GenIP}, with $G$ the solution map associated to the PDE \eqref{Eq:PDE} and $W$ a Gaussian white noise process indexed by $L^2(\Ocal)$.  An illustration of the problem with synthetic data is provided in Figure \ref{Fig:StatProbl} below. Among the numerous applications areas, inverse problems based on elliptic PDEs are important building blocks in oil reservoir modelling \cite{Y86}. Source identification problems have been extensively investigated in the applied mathematics and statistics communities; see \cite{adavani2010fast,elvetun2021regularization,GvdVY20,GK20} and the many references therein.

\begin{figure}[t]
\centering
\includegraphics[width=5.5cm,height=4.5cm]{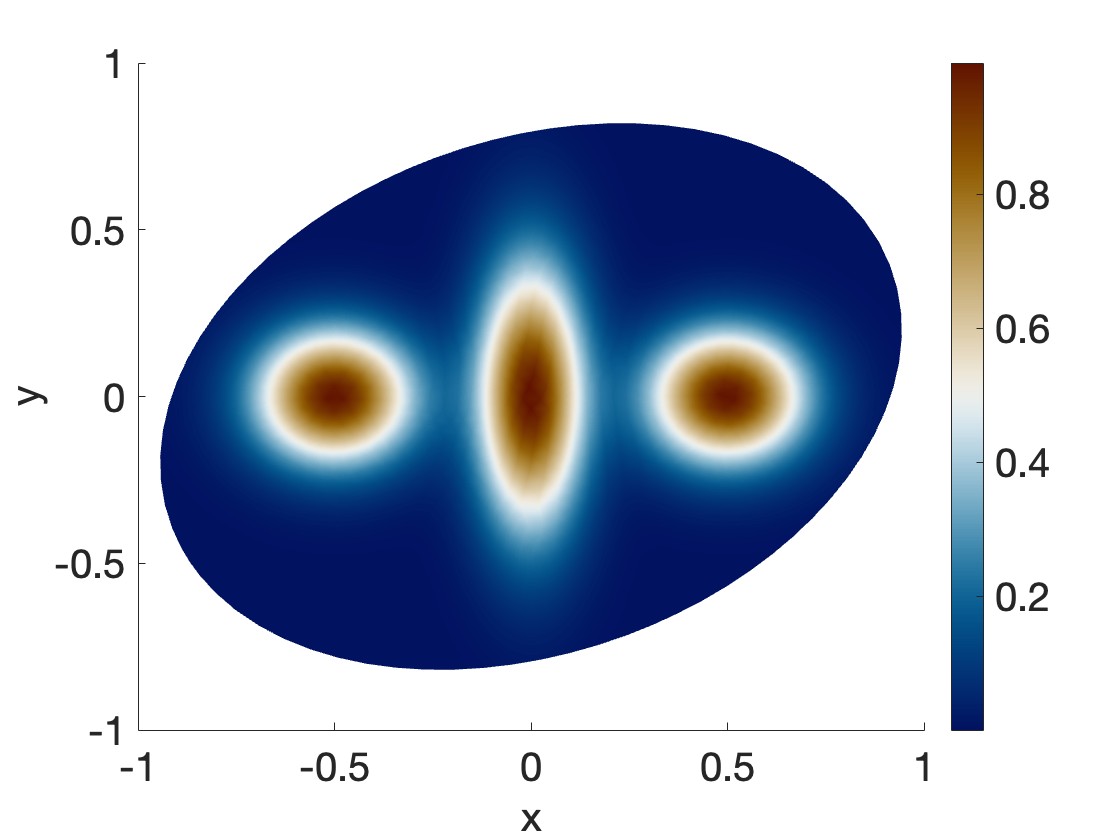}
\includegraphics[width=5.5cm,height=4.5cm]{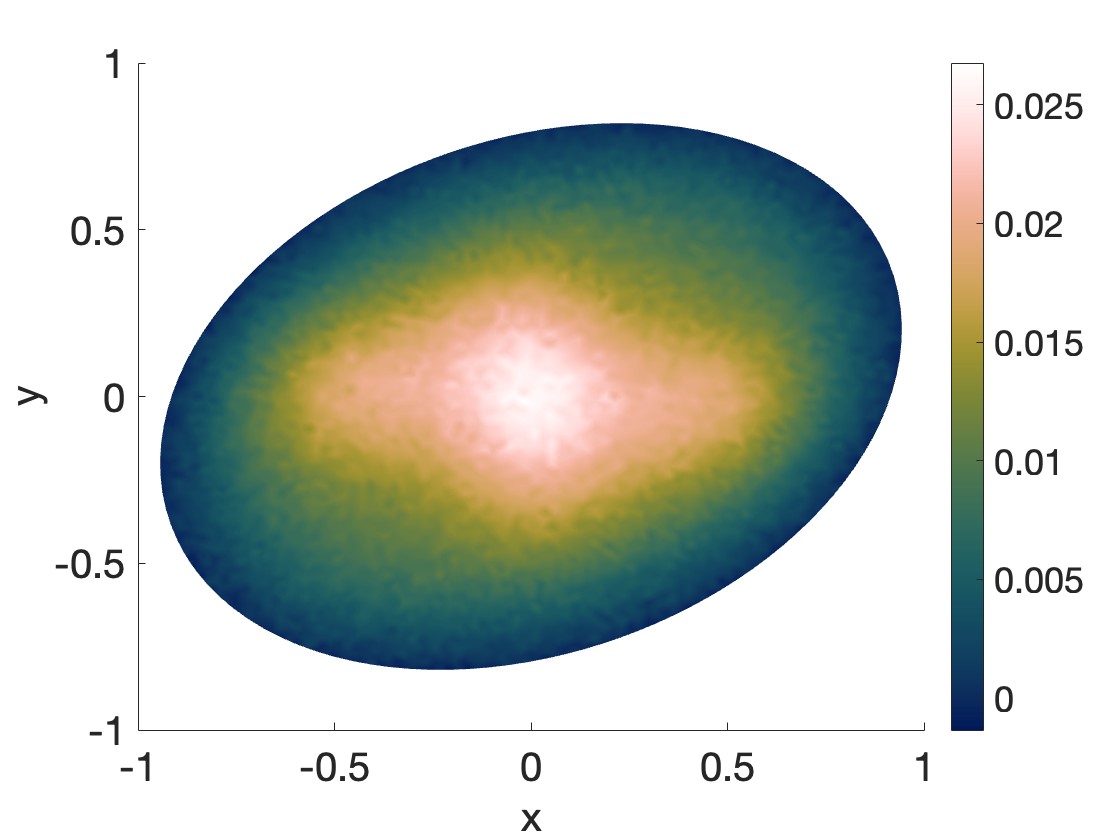}
\caption{Left: a source function $f$ on a rotated elliptically-shaped domain, modelling three heat sources centred at the points $(-0.5,0)$, $(0,0)$ and $(0,0.5)$. Right: noisy observations of the associated PDE solution $G(f)$ (with fixed diffusivity $c$).}
\label{Fig:StatProbl}
\end{figure}

\begin{figure}[t]
\centering
\includegraphics[width=5.5cm,height=4.5cm]{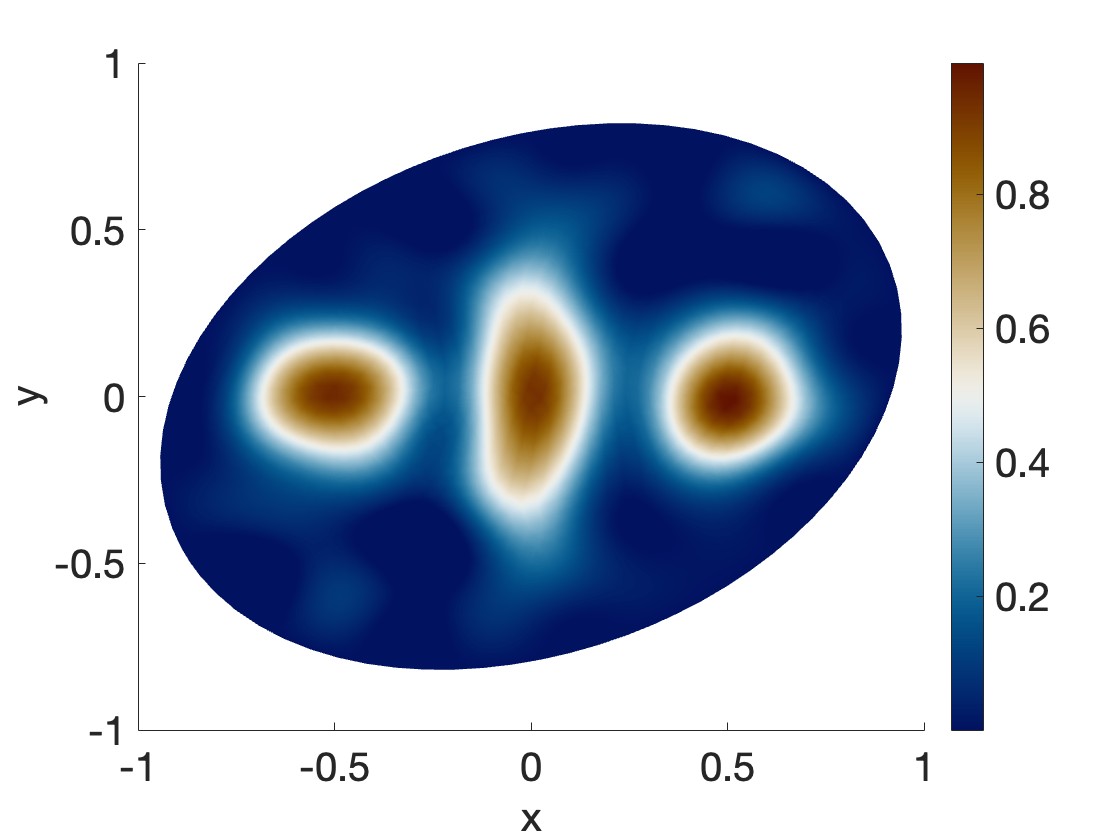}
\includegraphics[width=5.5cm,height=4.5cm]{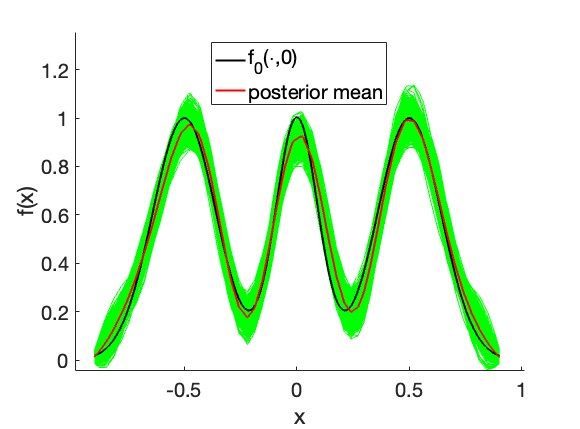}
\caption{Left: the posterior mean estimate of the source function $f$. Right: cross section along the $x$-axis of the true source function $f_0$ (in black), of the posterior mean (in red), and of $2500$ posterior samples (in green).}
\label{Fig:PostInf}
\end{figure}

	Here, we shall pursue the popular nonparametric Bayesian approach to inverse problems \cite{S10}. We shall assign to the unknown source $f$ a prior distribution $\Pi(\cdot)$ supported on the function space $L^2(\Ocal)$ and then, following the Bayesian paradigm, combine it with the (Gaussian) data likelihood induced by the statistical model \eqref{Eq:GenIP} to form the posterior distribution $\Pi(\cdot|Y^\eps)$ of $f|Y^\eps$, representing our updated belief about the inferential target and providing us with point estimates and uncertainty quantification. In particular, Gaussian priors represent a natural choice for inference in the context of the observation scheme \eqref{Eq:GenIP} due to their conjugacy, which leads to convenient explicit formulae for the posterior distribution.

	The Bayesian approach to inverse problems dates back at least to the 1980's, with early seminal work laid out in \cite{tarantola1982inverse,lehtinen1988statistical,lehtinen1989linear} among the others, and has since gained enormous popularity across several applied fields. See the monographs \cite{KS04,T05}, as well as the more recent reviews \cite{BB18,AMOS19}, where further references can be found. Over the last decade, a large number of articles have investigated the theoretical recovery performance of Gaussian priors in linear inverse problems. We refer in particular to \cite{KvdVvZ11,KvdVvZ13,ASZ14,KLS16,giordano2022nonparametric}, and also mention \cite{AN19,NvdGW20,GN20,MNP21a,giordano2023bayesian} for results in nonlinear problems. Recently, Giordano and Kekkonen \cite{GK20}, building on earlier findings by \cite{MNP19}, identified general conditions on the forward operator, on the ground truth and on the prior distribution under which `semiparametric' Bernstein-von Mises theorems can be obtained, characterising the asymptotic shape of the posterior distribution of a large collection of linear functionals of the unknown. For the problem of estimating the source $f$ in the elliptic PDE \eqref{Eq:PDE} from observations $Y^\eps$ arising as in \eqref{Eq:GenIP}, they showed that a wide class of `standard' centred Gaussian process priors (such as the one associated to the commonly used Matérn kernel) yield, in the small noise limit $\eps\to0$ and under the frequentist assumption that the data have been generated by a fixed ground truth $f_0$, valid and optimal estimation and uncertainty quantification, via the posterior mean estimate and credible intervals centred around it, for one-dimensional aspects of the unknown source function.

	For this study, we shall focus on the implementation of the Bayesian procedures with Gaussian priors for the inverse problem of source identification that were investigated in \cite{GK20}, corroborating the theory developed therein with a numerical simulation study. To this end, we will first review, in Section \ref{Sec:Theory}, the asymptotic results derived in \cite{GK20}. We will then provide, in Section \ref{Sec:Examples}, two examples of Gaussian priors satisfying the assumptions of the general theory; in particular, we will consider centred Gaussian series priors defined on the eigenbasis of the Dirichlet-Laplacian, as well as centred stationary Gaussian process priors associated to the Matérn covariance kernel. Both classes of priors are of practical interest and widely used, and they will lead to two different discretisation strategies for the implementation.

	We will present the numerical simulation study in Section \ref{Sec:Numerics}, were we will provide an implementation of posterior inference in the source identification problem for Gaussian series priors and for the Matérn process priors. The series approach will hinge on a natural discretisation of the parameter space via high-dimensional basis expansions, while for the Matérn priors we will employ  piecewise linear functions defined on the elements of a deterministic triangular mesh.  Under a suitable discretisation of the statistical model \eqref{Eq:GenIP}, we will provide the explicit formulae for the Gaussian conjugate posterior distributions, which we will exploit to efficiently compute the posterior mean estimates and to efficiently implement, via posterior sampling, credible sets for uncertainty quantification. In the study, we will numerically investigate the asymptotic concentration of the posterior distribution around the ground truth generating the data, the efficiency of the posterior mean estimator (in terms of minimality of its asymptotic variance), and further the frequentist coverage of the obtained credible intervals. Overall, the study shows a close correspondence between the numerical results and the expected performance predicted by the theory developed in \cite{GK20}. The posterior mean estimate (relative to a Gaussian series prior) and the associated uncertainty quantification (obtained via the cross-section of 2500 posterior samples along the $x$-axis) are depicted in Figure \ref{Fig:PostInf}, to be compared to the true source function shown in Figure \ref{Fig:StatProbl}. The reproducible (MATLAB) code used for the study is available at: \url{https://github.com/MattGiord/Bayesian-Source-Identification}.

%
%
%
%
%

\section{A nonparametric Bayesian approach with Gaussian priors}
\label{Sec:BayesApproach}

In this section, we formally introduce the statistical model and the Bayesian procedures of interest for the present article. We will provide details on the likelihood, the prior and the posterior distribution in Section \ref{Sec:PriorPost}. We will then review the asymptotic results of \cite{GK20} in Section \ref{Sec:Theory}, and further present two classes of Gaussian priors to which the theory applies in Section \ref{Sec:Examples}.

%
%
%

\subsection{Details on the Bayesian model}
\label{Sec:PriorPost}

Throughout, let $\Ocal\subset\R^d$, $d\in\N$, be a non-empty, open and bounded set with smooth boundary $\partial\Ocal$. For a fixed, smooth and positive `diffusion coefficient' $c\in C^\infty(\overline\Ocal)$, $\inf_{x\in\Ocal} c(x)>0$, let $G:L^2(\Ocal)\to L^2(\Ocal)$ be the (linear) source-to-solution map associated to the elliptic PDE \eqref{Eq:PDE}; see Appendix B in \cite{GK20} for details. For such $G$, and a given noise level $\eps>0$,
consider observations $Y^\eps$ arising as in \eqref{Eq:GenIP} for some unknown $f\in L^2(\Ocal)$, where $W$ is a white noise process indexed by $L^2(\Ocal)$ defined on some probability space $(\Omega,\Scal,\Pr)$, that is the centred Gaussian process $(W(g),\ g\in L^2(\Ocal))$ with covariance $\Enorm[W(g_1)W(g_2)]=\langle g_1,g_2\rangle_2$. Throughout most of the article, we will regard $\eps$ as known. In practice, it can often be replaced by an estimate (cf.~Section \ref{Sec:VarEstim}). As described, for example, in Chapter 1 of \cite{GN16}, observing $Y^\eps$ is understood as observing a realisation of the Gaussian process $(Y^\eps(g), \ g\in L^2(\Ocal))$ with mean $\Enorm[Y^\eps(g)]=\langle G(f),g\rangle_2$ and covariance $\Cov[Y^\eps(g_1),Y^\eps(g_2)]=\langle g_1,g_2\rangle_2$. Such observation scheme serves as a convenient continuous counterpart of the inverse regression model
\begin{equation}
\label{Eq:DiscrModel}
	Y_i = G(f)(X_i)+\sigma W_i, \qquad i=1,\dots,n,
\end{equation}
comprising noisy point evaluations of the PDE solution $G(f)$ over a set of points $X_1,\dots,X_n\in\Ocal$, corrupted by Gaussian measurement errors $\sigma W_1,\dots,\sigma W_n\iid N(0,\sigma^2)$, $\sigma>0$, known to be asymptotically equivalent (in the sense of \cite{BL96}) to the white noise model \eqref{Eq:GenIP} under suitable assumptions on the grid and the calibration $n/\sigma^2\simeq \eps^{-2}$.

	For continuous observations $Y^\eps$ arising as in \eqref{Eq:GenIP}, for any $f\in L^2(\Ocal)$, the (cylindrically-defined) law $P^\eps_f$ of $Y^\eps$ is absolutely continuous with respect to the law $P_0^\eps$ of the scaled white noise $\eps W$, with log-likelihood
$$ 
	\ell_\eps(f)
	:=\log \frac{dP^\eps_f}{dP^\eps_0}(Y^\eps) 
	= \frac{1}{\eps^2}Y^{\eps}[G(f)]-\frac{1}{2\eps^2}\|G(f)\|_2^2.  
$$
In view of the joint measurability of $\ell_\eps$, regarding $f$ as a random function with values in $L^2(\Ocal)$ and assigning to it any prior distribution in the form of a Borel probability measure $\Pi(\cdot)$ supported on $L^2(\Ocal)$ then induces, via Bayes' formula (for example, \cite[p.7]{GvdV17}), the posterior distribution
$$
	\Pi(A|Y^{\eps}) 
	= \frac{\int_A e^{\ell_\eps(f)}d\Pi(f)}
	{\int_{L^2(\Ocal)} e^{\ell_\eps(f')}d\Pi(f')},
	\qquad A\subseteq L^2(\Ocal)\ \textnormal{measurable},
$$ 
that is, the conditional distribution of $f|Y^\varepsilon$. In particular, it will be of interest to consider Gaussian priors which, in view of the linearity of the forward operator $G$ and the normal assumption on the noise $W$, will lead to conjugate Gaussian posteriors. Concrete formulae are provided in Section \ref{Sec:Numerics}. In the following, we will repeatedly use elements of the theory of Gaussian processes and measures on Hilbert spaces, and we refer to \cite[Chapter 2]{GN16} for the necessary background. For a Gaussian prior $\Pi(\cdot)$ on $L^2(\Ocal)$, the `information geometry' is encoded within an associated reproducing kernel Hilbert space (RKHS) of functions defined on the domain $\Ocal$, strictly contained inside $L^2(\Ocal)$. Popular prior choices in applications and theoretical studies typically model functions belonging to a `smoothness scale', with associated RKHS equal to (or included in) a Sobolev space $H^\alpha(\Ocal)$, for some regularity level $\alpha>0$. These include Gaussian series priors defined on bases spanning the Sobolev scale, as well as stationary Gaussian processes with the Matérn covariance kernel; see \cite[Chapter 11]{GvdV17}. 

%
%
%

\subsection{Theoretical guarantees for estimation and uncertainty quantification}
\label{Sec:Theory}

The asymptotic properties of nonparametric Bayesian procedures with Gaussian priors in inverse problems have recently been investigated by Giordano and Kekkonen \cite{GK20}, resulting, under general `regularity conditions' for the forward operator, for the ground truth and for the prior distribution, in semiparametric Bernstein-von Mises theorems that entail the convergence of certain one-dimensional posterior distributions to limiting Gaussian probability measures with minimal variance, in the small noise limit and under the frequentist assumption that the data have been generated by a fixed ground truth. These results were then leveraged in \cite{GK20} to prove the asymptotic efficiency of the associated posterior mean estimators, as well as to derive theoretical guarantees certifying that credible intervals centred around them are asymptotically valid confidence intervals with minimal width. In this section, we provide a review of the findings of \cite{GK20} for the inverse problem of source identification. These will later be corroborated by the results of the numerical simulation study presented in Section \ref{Sec:Numerics}.

	The investigation of \cite{GK20} builds on the semiparametric approach to the Bernstein-von Mises theorem in infinite-dimensional statistical models developed by Castillo and Nickl \cite{CN13} and later refined by Monard et al.~\cite{MNP19} in the inverse problem setting. It is based on the study of the posterior distributions of a class of  scaled and centred one-dimensional functionals of the unknown which, in the context of the source identification problem, take the form $\eps^{-1}\langle f- \bar{f}_\eps,\psi\rangle_2$ for test functions $\psi\in L^2(\Ocal)$, where $\bar{f}_\eps:=E^\Pi[f|Y^\eps]$ is the posterior mean. Let $\Lcal(\eps^{-1}\langle f- \bar{f}_\eps ,\psi\rangle_2|Y^{\eps})$ denote the associated scaled and centred posterior distribution.

\begin{theorem}[Theorem 4.1 in \cite{GK20}]\label{Theo:BvM}
Let $\Pi(\cdot)$ be a centred Gaussian Borel probability measure supported on $L^2(\Ocal)$ with RKHS equal to $H^\alpha(\Ocal)$ for some $\alpha>d/2$. Let $f_0\in H^\beta(\Ocal)$, for some $\beta>\alpha-d/2$, be compactly supported inside $\Ocal$, and consider observations $Y^\eps\sim P^\eps_{f_0}$ from the statistical model \eqref{Eq:GenIP} with $G(f)$ the solution to the PDE \eqref{Eq:PDE} and $f=f_0$. Then, for any $\gamma>2+d/2$ and any compactly supported test function $\psi\in H^\gamma(\Ocal)$, we have
\begin{equation}
\label{Eq:BvM}
	\Lcal(\eps^{-1}\langle f- \bar{f}_\eps ,\psi\rangle_2|Y^{\eps})
	\overset{\Lcal}{\longrightarrow} N(0,\|\nabla\cdot(c\nabla \psi)\|^2_2),
\end{equation}
in $P^\eps_{f_0}$-probability as $\eps\to0$.
\end{theorem}

	The result asserts that the random (data-dependent) one-dimensional probability distribution $\Lcal(\eps^{-1}\langle f- \bar{f}_\eps ,\psi\rangle_2|Y^{\eps})$ converges (in the topology of weak convergence) in probability to a centred normal distribution with variance $\|\nabla\cdot(c\nabla \psi)\|^2_2$. The latter can be shown to be minimal \cite[Remark 2.4]{GK20}, as it coincides with the Cramér-Rao lower bound for estimating the one-dimensional quantity $\langle f,\psi\rangle_2$ from data $Y^\varepsilon$ arising as in \eqref{Eq:GenIP}. Furthermore, the class of test functions $\psi$ for which the convergence \eqref{Eq:BvM} is obtained is to be understood to be maximal, in the sense that the infinite-dimensional Gaussian probability measure with marginals identified by the right hand side of \eqref{Eq:BvM} is tight (a necessary condition for weak convergence) only when $\gamma> 2+d/2$; see Lemma 4.2 in \cite{GK20} and the related discussion.

 	A first important consequence of Theorem \ref{Theo:BvM} is a central limit for the `plug-in' posterior mean estimators $\langle \bar f_\eps,\psi\rangle_2$ of the one dimensional aspects $\langle f,\psi\rangle_2$ of the unknown. Note that, for Gaussian priors, these can be efficiently computed via the explicit formulae for the conjugate Gaussian posteriors. The central limit follows, as argued in Remark 2.4 in \cite{GK20}, from the convergence of moments in the limit \eqref{Eq:BvM}. In particular, under the assumptions of Theorem \ref{Theo:BvM}, it holds that
\begin{equation}
\label{Eq:CLT}
	\eps^{-1}\left(\langle \bar{f}_\eps ,\psi\rangle_2 
	- \langle f_0,\psi\rangle_2\right)
	 \overset{d}{\longrightarrow} 
	N(0,\|\nabla\cdot(c\nabla\psi)\|^2_2),
\end{equation}
as $\eps\to0$. In view of the aforementioned minimality of the asymptotic variance $\|\nabla\cdot(c\nabla\psi)\|^2_2$, the result indeed implies the asymptotic efficiency of the plug-in estimators $\langle \bar f_\eps,\psi\rangle_2$.

	The second implication of the Berstein-von Mises result stated in Theorem \ref{Theo:BvM} concerns the coverage and width of credible intervals built around the efficient estimators $\langle \bar f_\eps,\psi\rangle_2$, which can be shown to be asymptotically valid frequentist confidence intervals and to have diameter shrinking at the optimal parametric rate $\eps^{-1}$. For any level $a\in(0,1)$, consider the $(1-a)\%$-credible interval
\begin{equation}
\label{Eq:CredInt}
	C_{\eps,a}
	:= \{z\in\R:  | z -\langle \bar{f}_\eps ,\psi\rangle_2|
	\le R_{\eps,a}\},
\end{equation}
where $R_{\eps,a}>0$ is the $(1-a/2)\%$-quantile of the one-dimensional (Gaussian) posterior distribution of $\langle f ,\psi\rangle_2|Y^\eps$, so that
$$
	\Pi\left(f:\langle f,\psi\rangle_2
	\in C_{\eps,a}|Y^{\eps}
	\right) 
	= 1-a.
$$
Then, in the setting of Theorem \ref{Theo:BvM}, the asymptotic frequentist coverage of $C_{\eps,a}$ is given by
\begin{equation}
\label{Eq:Coverage}
	P^\eps_{f_0}
	\left( \langle f_0,\psi\rangle_2\in C_{\eps,a}\right) \to 1-a,
\end{equation}
as $\eps\to0$, while its radius $R_{\eps,a}$ satisfies
$$
 	R_{\eps,a} = O_{P^\eps_{f_0}}(\eps^{-1}).
$$
See Corollary 2.5 in \cite{GK20}. Note that although an analytic formulation of the credible intervals $C_{\eps,a}$ requires the derivation of the quantiles of the one dimensional posterior distributions $\langle f,\psi\rangle_2|Y^\eps$, these can typically be numerically approximated by efficiently sampling from the explicitly available conjugate posterior distributions.
	
%
%
%

\subsection{Examples of Gaussian priors}
\label{Sec:Examples}

In this section, we provide two concrete examples of Gaussian priors to which Theorem \ref{Theo:BvM} applies. For both instances, an implementation of the resulting posterior inference will be presented in Section \ref{Sec:Numerics} below, based on two different discretisation strategies. The first example concerns Gaussian series priors.

\begin{example}[Gaussian series priors on the Dirichlet-Laplacian eigenbasis]\label{Ex:DirichLapl}
Let $(\phi_j, \ j\in\N)\subset H^1(\Ocal)\cap C^\infty(\overline\Ocal)$ be the orthonormal basis of the space $L^2(\Ocal)$ formed by the eigenfunctions of the (negative) Dirichlet-Laplacian,
\begin{equation}
\label{Eq:Eigen}
\begin{split}
-\Delta \phi_j -\lambda_j \phi_j &=0, \ \ \textnormal{on}\ \ \Ocal \\
\phi_j&=0, \ \ \textnormal{on}\ \ \partial\Ocal,
\end{split}
\qquad j\in\N,
\end{equation}
with associated eigenvalues $0<\lambda_1<\lambda_2\le\lambda_3\le\dots, $ satisfying $\lambda_j\to\infty$ as $j\to\infty$ according to Weyl's asymptotics, namely $\lambda_j=O(j^{2/d})$ as $j\to\infty$. We refer to Example 6.3 and Section 7.4 in \cite{HT07} for details. The associated Hilbert scale
$$
	\H^\alpha
	:=\Bigg\{f\in L^2(\Ocal) : \|f\|^2_{\H^\alpha}:=\sum_{j=1}^\infty \lambda_j^\alpha
	|\langle f,\phi_j\rangle_2|^2
	<\infty\Bigg\},
	\qquad \alpha\ge0,
$$
then satisfies $\H^0=L^2(\Ocal)$ (with equality of norms) and the continuous (strict) embedding $\H^\alpha\subset H^\alpha(\Ocal)$ for all $\alpha>0$ \cite[p.~472]{T11}. In fact, it holds that $\|f\|_{\H^\alpha}\simeq \|f\|_{H^\alpha}$ for all $f\in \H^\alpha$ and $\alpha\ge0$.

	For any $\alpha>d/2$, consider the Gaussian random series
\begin{equation}
\label{Eq:SeriesPrior}
		F:= \sum_{j=1}^\infty \lambda_j^{-\alpha/2} F_j \phi_j,
		 \qquad F_j\iid N(0,1),
\end{equation}
corresponding to the Karhunen-Loève expansions of certain commonly used Gaussian process priors with covariance kernel given by an inverse power of the Laplacian \cite[Section 2.4]{S10}. By Weyl's asymptotics, we have
$$
	\Enorm[\| F\|_2^2] =  \sum_{j=1}^\infty \lambda_j^{-\alpha} 
	\simeq \sum_{j=1}^\infty j^{-2\alpha/d}
	<\infty,
$$
since $2\alpha/d>1$, showing that $F$ takes values almost surely in $L^2(\Ocal)$. By Lemma I.5 in \cite{GvdV17}, the law $\Pi(\cdot)$ of $F$ is then seen to define a Gaussian Borel probability measure supported on $L^2(\Ocal)$. Furthermore, by Theorem I.12 in \cite{GvdV17}, its RKHS is equal to $\H^\alpha$. Noting that, for any compactly supported test function $\psi\in H^\gamma(\Ocal)$, $\gamma>2 + d/2$, the approximation argument in the proof of Theorem 4.1 in \cite{GK20} can be carried out with minimal modifications with $H^\alpha(\Ocal)$ replaced by $\H^\alpha$, we conclude that Theorem \ref{Theo:BvM} applies by modelling the unknown source function $f$ via the Gaussian series \eqref{Eq:SeriesPrior}.

	While not explicitly available for general domains $\Ocal$, we note that the Dirichlet-Laplacian eigenbasis can be numerically computed via efficient finite element methods for elliptic eigenvalue problems, offering a broadly applicable framework for implementation. More details will be provided in Section \ref{Sec:Numerics} below.
\end{example}

	A second example of interest involves stationary Gaussian processes defined via a  covariance kernel of choice. In particular, the Matérn kernel is widely used in applications \cite[Section 4.2]{RW06}.

\begin{example}[Matérn process priors]\label{Ex:Matern}
	Let $F=(F(x),\ x\in\Ocal)$ be the centred and stationary Gaussian process with Matérn covariance kernel
\begin{equation}
\label{Eq:MatCov}
	C_{\alpha,\ell}(x,y)= \frac{2^{1-\alpha}}{\Gamma(\alpha)}\left(\frac{|x-y|\sqrt {2 \alpha}}{\ell}\right)^\alpha
	B_\alpha \left(\frac{|x-y|\sqrt{2\alpha}}{\ell}\right),
	\qquad x,y\in\Ocal,
\end{equation}
with regularity parameter $\alpha>d/2$ and length scale $\ell>0$. Above, $\Gamma$ denotes the gamma function and $B_\alpha$ is the modified Bessel function of the second kind. The finite dimensional distributions of $F$ are identified by the relation
\begin{equation}
\label{Eq:FDDs}
	(F(x_1),\dots,F(x_M))^T\sim N_M(0,\mathbf C),
\end{equation}
with $\mathbf C:=(C_{\alpha,\ell}(x_h,x_m))_{h,m=1}^M\in \R^{M,M}$, holding for any $M\in\N$ and any $x_1,\dots, x_M\in\Ocal$. By Lemma I.4 in \cite{GvdV17}, a version of $F$ can be identified with sample paths belonging almost surely to the H\"older space $C^{\alpha'}(\Ocal)\subset L^2(\Ocal)$ for any $0<\alpha'<\alpha-d/2$, and therefore, in view of Lemma I.7 in \cite{GvdV17}, the law $\Pi(\cdot)$ of such version defines a Gaussian Borel probability measure supported on $L^2(\Ocal)$. Moreover, the results in Section 11.4.4 in \cite{GvdV17} imply that the RKHS of $F$ equals, with norm equivalence, the set of restrictions to the domain $\Ocal$ of functions in the Sobolev space $H^\alpha(\R^d)$. Since $\Ocal$ is assumed to have smooth boundary, the latter is indeed equal to $H^\alpha(\Ocal)$. Thus, Theorem \ref{Theo:BvM} applies with $\Pi(\cdot)$ a Matérn process prior with covariance kernel \eqref{Eq:MatCov}.
\end{example}

%
%
%
%
%

\section{Numerical simulation study}
\label{Sec:Numerics}

For illustration, we take as working domain the area $\Ocal$ contained inside a rotated ellipse with horizontal semi-axis of unit length, vertical semi-axis of length $3/4$, and rotation angle $\theta=\pi/6$,
$$
	\{(\cos(t)\cos(\theta) - 3/4\sin(t)\sin(\theta),3/4\sin(t)\cos(\theta)+\cos(t)\sin(\theta)),\
	t\in[0,2\pi)\}.
$$
For an unknown source function $f\in L^2(\Ocal)$, we assume in practice that we are given $n$ noisy point evaluations $\mathbf Y:=(Y_1,\dots,Y_n)^T\in \R^n$ of the solution $G(f)$ to the PDE \eqref{Eq:PDE} generated according to the equivalent discrete statistical model \eqref{Eq:DiscrModel}, for a given deterministic grid of points $x_1,\dots,x_n\in\Ocal$ comprising the nodes of a triangular mesh covering $\Ocal$ (Figure \ref{Fig:Eigenfun}, top-left). We then seek to estimate $f$ from data $\mathbf Y$.

%
%
%


\subsection{Posterior inference with Gaussian series priors}
\label{Sec:SeriesPriors}

%

\subsubsection{Methodology}
\label{Sec:Methodology}

For the Gaussian series priors defined via the Dirichlet-Laplacian eigenpairs $\{(\phi_j,\lambda_j),\ j\in\N\}$ considered in Example \ref{Ex:DirichLapl}, we discretise the parameter space by modelling the unknown source function $f$ as the finite sum
\begin{equation}
\label{Eq:Discretisation}
	f = \sum_{j=1}^J f_j \phi_j, \qquad f_1,\dots,f_J \in\R,\qquad J\in\N.
\end{equation}
For any such $f$, the linearity of the forward map $G$ then implies that the discrete observations are given by
\begin{align*}
	Y_i
	=\sum_{j=1}^J f_j G(\phi_j)(x_i)+\sigma W_i
	\equiv (\mathbf G \mathbf f)_i +\sigma W_i
\end{align*}
where $\mathbf G:=[G(\phi_j)(x_i),\ i=1,\dots,n,\ j=1,\dots,J]\in \R^{n,J}$ and $\mathbf f:=(f_1,\dots,f_J)^T\in \R^J$, whereby the inverse regression model \eqref{Eq:DiscrModel} can be written in matrix notation as
\begin{equation}
\label{Eq:MatrixModel}
	\mathbf Y = \mathbf G\mathbf f +\sigma \mathbf W
\end{equation}
with $\mathbf W:=(W_1,\dots,W_n)^T\sim N_n(0,\mathbf I_n)$. Thus, for any given $\mathbf f\in\R^J$, $\mathbf Y|\mathbf f\sim N_n(\mathbf G\mathbf f,\sigma^2 \mathbf I_n)$.

\begin{figure}[t]
\centering
\includegraphics[width=5.5cm,height=4.5cm]{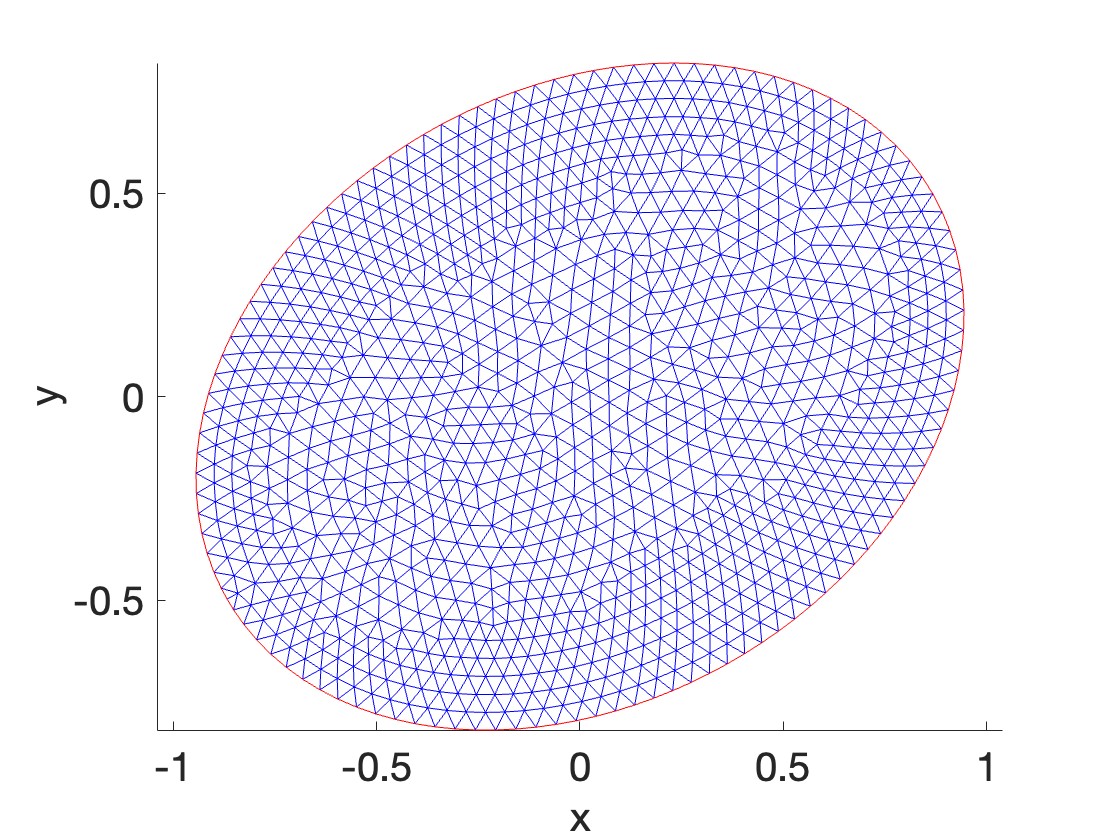}
\includegraphics[width=5.5cm,height=4.5cm]{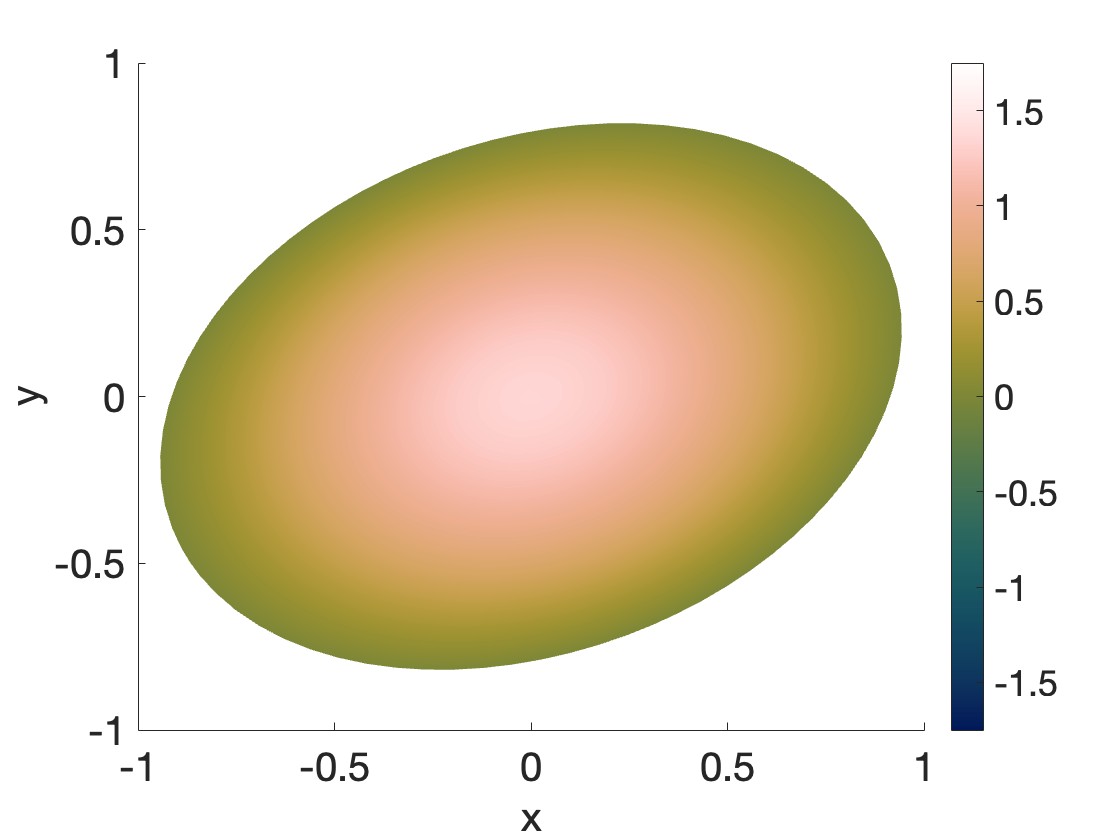}
\includegraphics[width=5.5cm,height=4.5cm]{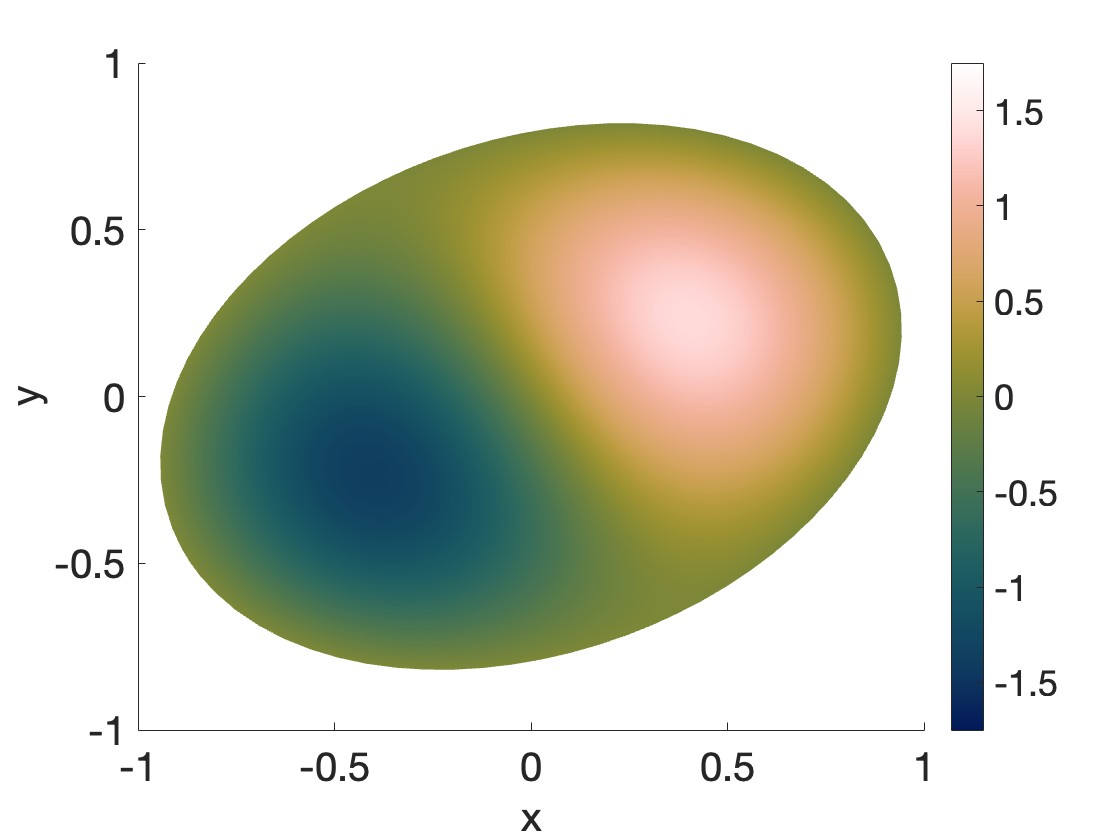}
\includegraphics[width=5.5cm,height=4.5cm]{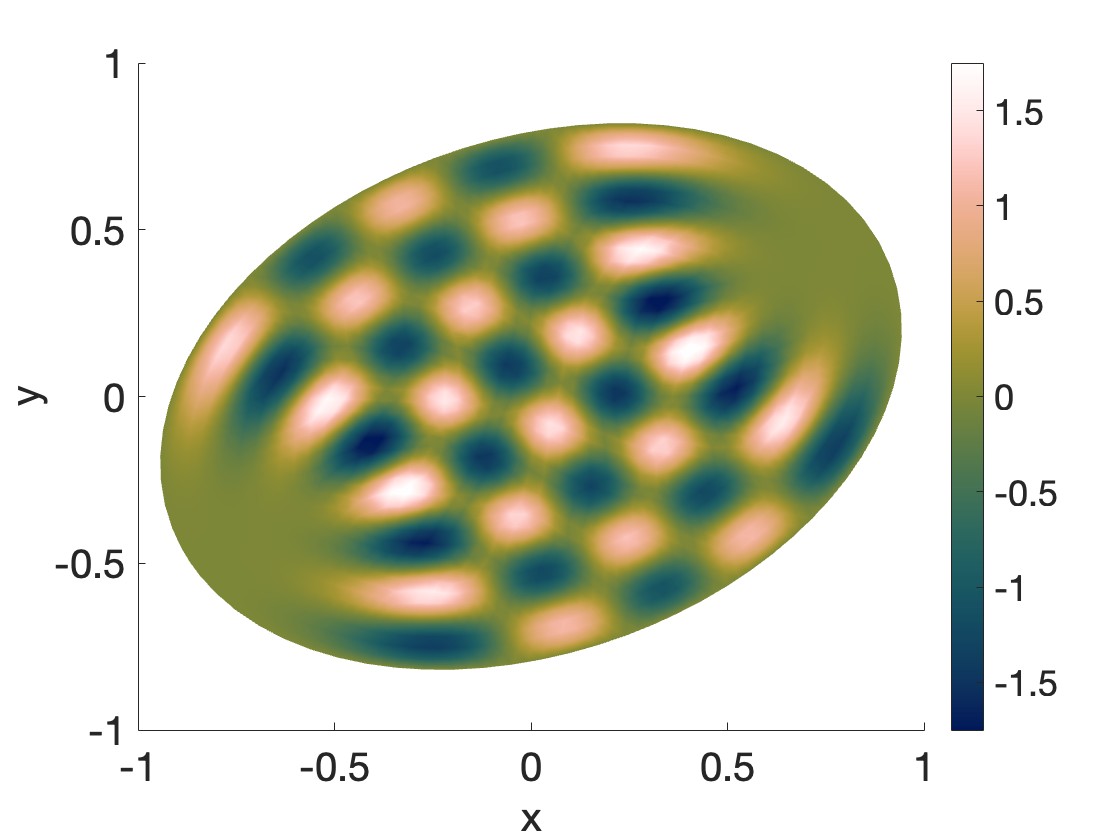}
\caption{Top-left: triangular mesh on a rotated elliptically-shaped domain with $n=4500$ nodes. Top-right, bottom-left and bottom-right, respectively: numerical approximations of the first Dirichlet-Laplacian eigenfunction $\phi_1$, of the second eigenfunction $\phi_2$ and of the eigenfunction $\phi_{84}$ associated to the largest eigenvalue $\lambda_{84}=493.3725$ found in the range $[0,500]$.}
\label{Fig:Eigenfun}
\end{figure}

\begin{figure}[t!h]
\centering
\includegraphics[width=5.5cm,height=4.5cm]{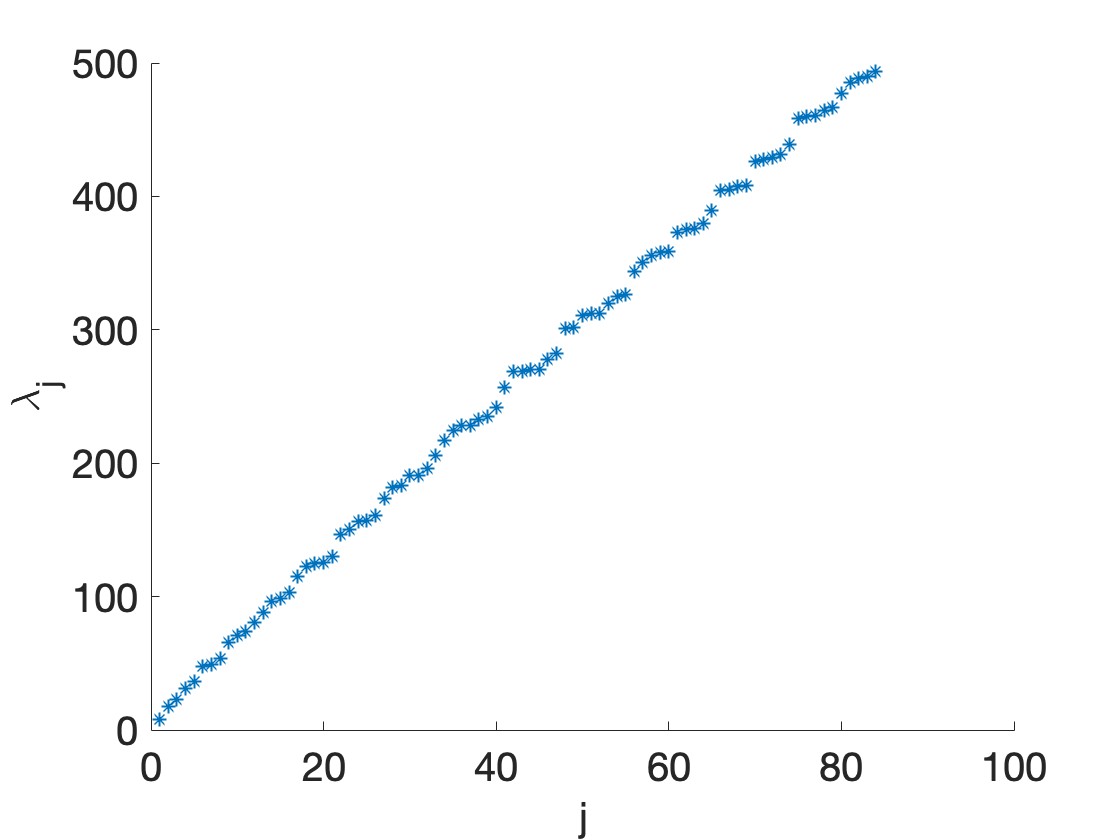}
\caption{Numerical approximations of the Dirichlet-Laplacian eigenvalues in the range $[0,500]$.}
\label{Fig:Eigenval}
\end{figure}

	In practice, outside certain special cases for the domain $\Ocal$ (such as squared or circular ones), the Dirichlet-Laplacian eigenpairs are not explicitly available. For general domains, we then resort to numerical methods to solve the elliptic eigenvalue problem \eqref{Eq:Eigen}. In particular, we used MATLAB PDE Toolbox, that allows to input a range of search $[0,\lambda_{\textnormal{max}}]$, $\lambda_{\textnormal{max}}>0$, for the eigenvalues, and then returns numerical approximations, obtained via finite element methods, of the eigenvalues in the prescribed interval and of the corresponding eigenfunctions. For the considered rotated elliptically-shaped domain, Figure \ref{Fig:Eigenfun} shows the first, the second and the last eigenfunction returned by the elliptic PDE solver. The range was set to $[0,\lambda_\textnormal{max}]=[0,500]$, for which $J=84$ eigenvalues were found. The (numerical approximations to) the eigenvalues are displayed in Figure \ref{Fig:Eigenval}. They exhibit a linear growth as expected from Weyl's asymptotics for bi-dimensional domains. The computation of the matrix $\mathbf G$ in the discretised observation model \eqref{Eq:MatrixModel} is performed by numerically solving, again using MATLAB PDE Toolbox, the elliptic PDE \eqref{Eq:PDE} with $f$ replaced by $\phi_j$ (or, more precisely, by the finite element approximation of $\phi_j$) and then evaluating $\mathbf G_{ij}:=G(\phi_j)(x_i)$ for all $i=1,\dots,n$ and $j=1,\dots,J$.

	Under the discretisation \eqref{Eq:Discretisation}, the Gaussian series prior described in Example \ref{Ex:DirichLapl} is approximately implemented by truncating the random series \eqref{Eq:SeriesPrior} at level $J$, and then assigning to the coefficients $f_1,\dots,f_J$ in \eqref{Eq:Discretisation} independent Gaussian priors $N(0,\lambda_j^{-\alpha})$, $j=1,\dots,J$. In the discretised observation model \eqref{Eq:MatrixModel}, this corresponds to assigning to the vector $\mathbf f$ the $J$-dimensional Gaussian prior with diagonal covariance matrix
\begin{equation}
\label{Eq:DiscrPrior}
	\mathbf f\sim N_J(0,\mathbf \Lambda), 
	\qquad \mathbf \Lambda:=\textnormal{diag}(\lambda_1^{-\alpha},
	\dots,\lambda_J^{-\alpha})\in\R^{J,J}.
\end{equation}
Thus, recalling that, according to \eqref{Eq:MatrixModel}, $\mathbf Y|\mathbf f\sim N_n(\mathbf G\mathbf f,\sigma^2 \mathbf I_n)$, a standard conjugate computation for multivariate models with Gaussian likelihood and prior yields the posterior distribution
\begin{equation}
\label{Eq:Conjugacy}
	\mathbf f|\mathbf Y\sim N_J(\bar{\mathbf f}_n,\mathbf \Lambda_n),
\end{equation}
where
\begin{equation}
\label{Eq:PostParam}
	\bar{\mathbf f}_n:=\frac{1}{\sigma^2}\mathbf \Lambda_n\mathbf G^T\mathbf Y; 
	\qquad \mathbf \Lambda_n:=(\sigma^{-2}\mathbf G^T\mathbf G+\mathbf \Lambda^{-1})^{-1}.
\end{equation}
Using the conjugate formulae \eqref{Eq:Conjugacy} and \eqref{Eq:PostParam}, it is straightforward to compute posterior mean estimates and drawing posterior samples. In turn, this allows to efficiently implement credible sets centred around the posterior mean, replacing the theoretical posterior quantiles (for example, the ones used in the definition of the credible intervals \eqref{Eq:CredInt}) with the empirical quantiles associated to a sufficiently large sample from the posterior distribution.

%

\subsubsection{Experiments}
\label{Sec:SeriesExp}

Throughout the numerical simulation study, the true source function (shown in Figure \ref{Fig:StatProbl}, left) was taken to be
$$
	f_0(x,y) = e^{-(5x-2.5)^2-(5y)^2} + e^{-(7.5x)^2-(2.5y)^2} + e^{-(5x-2.5)^2-(5y)^2},
	\qquad (x,y)\in\Ocal.
$$
Figure \ref{Fig:StatProbl} (right) shows $n=4500$ discrete noisy observations, over the nodes of the triangular mesh depicted in Figure \ref{Fig:Eigenfun} (top-left), of the corresponding PDE solution $G(f_0)$ arising as in the inverse regression model \eqref{Eq:DiscrModel} with noise standard deviation $\sigma=0.0005$ (with corresponding signal-to-noise ratio $\|G(f_0)\|_2/\sigma=37.55$).  The diffusion coefficient was taken to be $c(x,y):=2+5e^{-(5x-2)^2-(5y-2)^2}+5e^{-(5x+2)^2-(5y+2)^2}$, $(x,y)\in\Ocal$. The PDE solution $G(f_0)$ was calculated using MATLAB PDE Toolbox, which also contains the routine to create the triangular mesh.

\begin{figure}[!th]
\centering
\includegraphics[width=5.5cm,height=4.5cm]{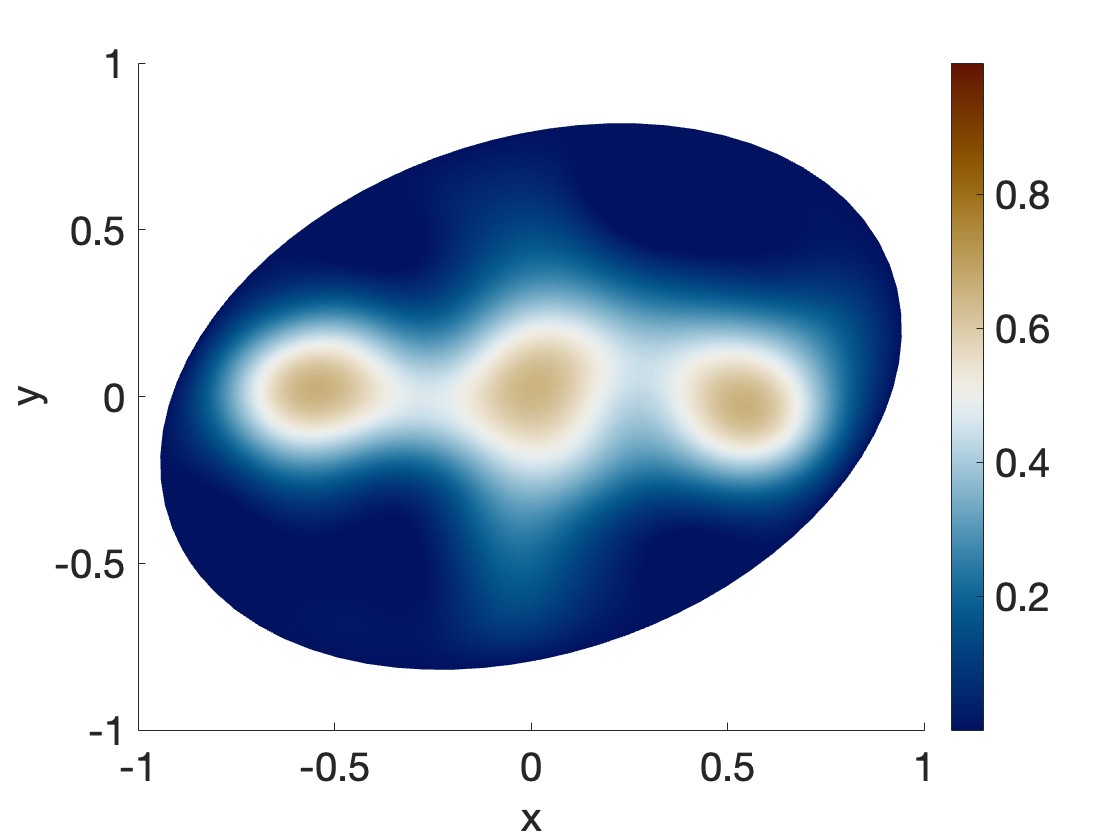}
\includegraphics[width=5.5cm,height=4.5cm]{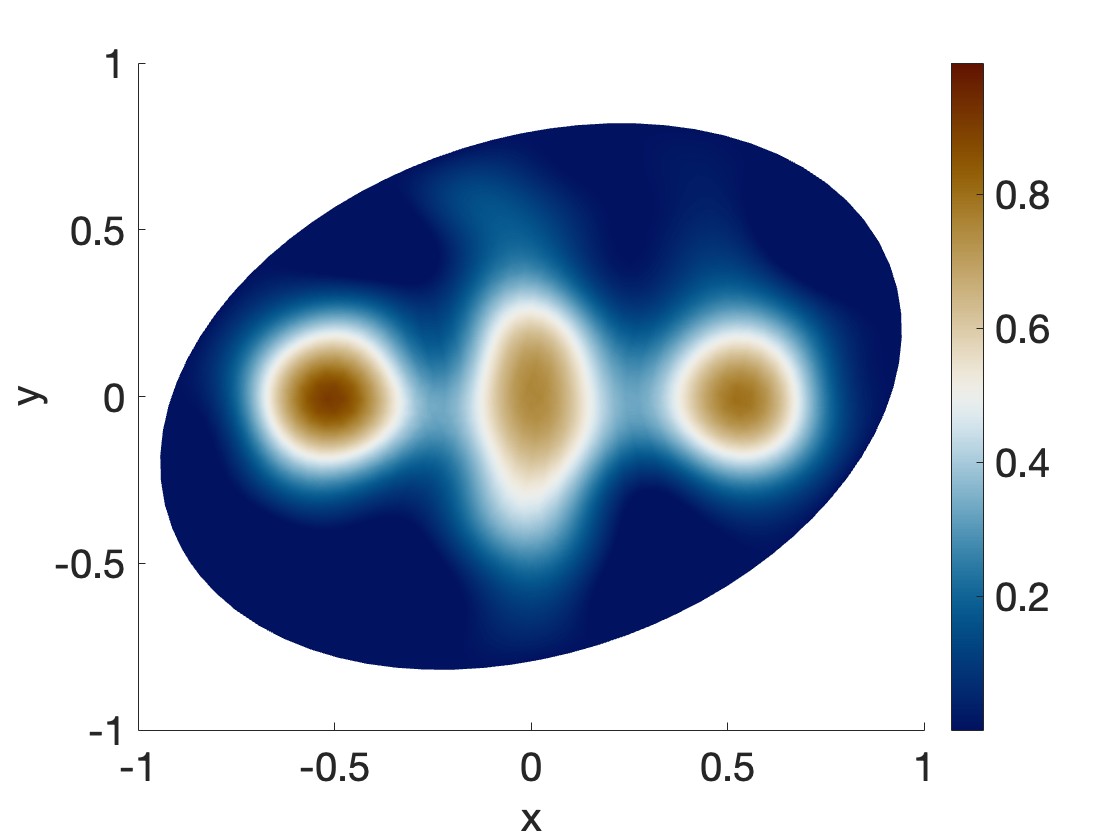}
\includegraphics[width=5.5cm,height=4.5cm]{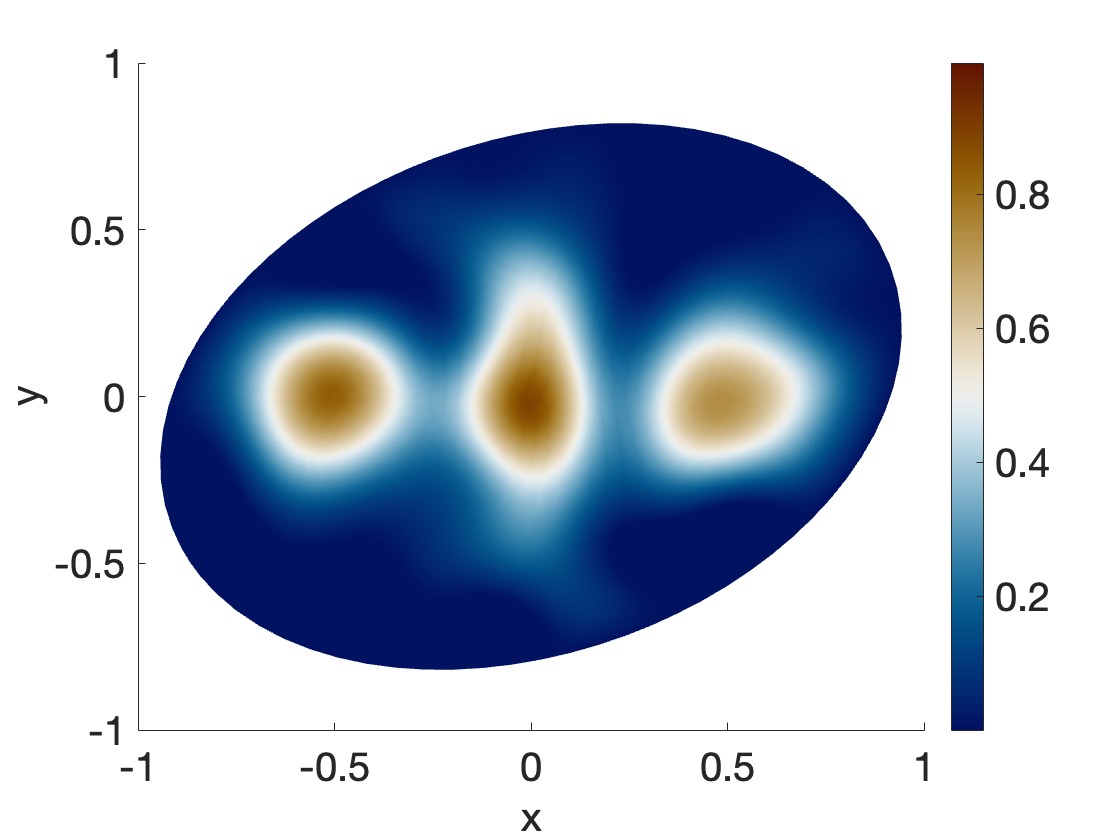}
\includegraphics[width=5.5cm,height=4.5cm]{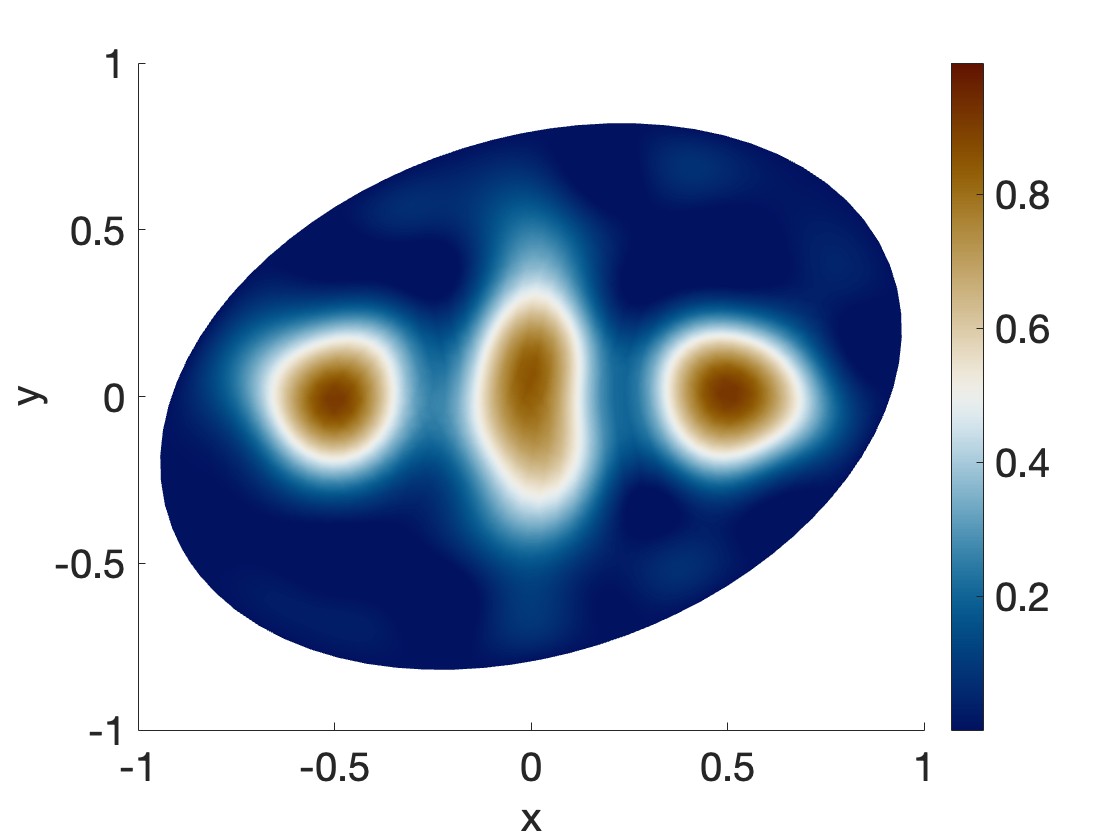}
\caption{Left to right, top to bottom: posterior mean estimates $\bar f_n$ of the source function $f$ for increasing sample sizes $n = 100,250,500,2000$.}
\label{Fig:Asymptotics}
\end{figure}

	The posterior mean estimate $\bar f_n:=\sum_{j=1}^J\bar{\mathbf f}_{n,j}\phi_j$ shown in Figure \ref{Fig:PostInf} (left) was obtained by computing the vector of coefficients $\bar{\mathbf f}_n$ according to the conjugate formula in \eqref{Eq:PostParam}. A diagonal Gaussian prior as in \eqref{Eq:DiscrPrior} was used, with regularity parameter $\alpha=3/4$. The parameter space was discretised using $J=84$ basis functions. The obtained $L^2$-estimation error was $\|\bar f_n - f_0\|_2=0.060077$. For comparison, $\|f_0\|_2=0.4764$ (with corresponding relative error $\|\bar f_n - f_0\|_2/\|f_0\|_2=12.5\%$), while the $L^2$-approximation error incurred by projecting $f_0$ onto the linear space spanned by the employed set of basis functions (furnishing a lower bound for the $L^2$-estimation error) is $0.0486$. The $2500$ posterior draws whose cross-sections along the $x$-axis are shown in Figure \ref{Fig:PostInf} (right) were sampled from the conjugate Gaussian posterior distribution in \eqref{Eq:Conjugacy}.

	Figure \ref{Fig:Asymptotics} provides an illustration of asymptotic convergence in the infinitely-informative data limit, showing the posterior mean estimates obtained for increasing sample sizes. The (decreasing) $L^2$-estimation errors for sample sizes ranging between $n=50$ and $n=4500$ are reported in Table \ref{Tab:Errors}. Across the experiments, the same discretisation with $J=84$ basis function and the same diagonal Gaussian prior with regularity $\alpha=3/4$ were used.

\begin{table}[!t]
\caption{$L^2$-estimation errors achieved by the posterior mean estimator $\bar f_n$ for increasing sample sizes.}
\label{Tab:Errors} 
\begin{tabular}{p{3cm}p{.925cm}p{.925cm}p{.925cm}p{.925cm}p{.925cm}p{.925cm}p{.925cm}p{.95cm}p{.7cm}}
\hline\noalign{\smallskip}
 $n$ & 50 & 100 & 250 & 500 & 750 & 1000 &2000 & 3000 & 4500  \\
\noalign{\smallskip}\hline\noalign{\smallskip}
$\|\bar f_n - f_0\|_2$ & 0.22  & 0.18& 0.13 & 0.099 & 0.088 & 0.078 & 0.076 & 0.069 & 0.060\\
\noalign{\smallskip}\hline\noalign{\smallskip}
$\|\bar f_n - f_0\|_2/\|f_0\|_2$ & 45.8\%  & 37.5\% & 27.1\%  & 20.6\% & 18.3\%&  16.3\% &15.8\%& 14.3\%& 12.5\% \\
\noalign{\smallskip}\hline\noalign{\smallskip}
\end{tabular}
\end{table}

	We next consider semiparametric inference for one-dimensional linear functionals $\langle f,\psi\rangle_2$, $\psi\in L^2(\Ocal)$, and provide a numerical illustration of the asymptotic results presented in Section \ref{Sec:Theory}. In particular, we focus on test functions $\psi=\phi_j$ , $j\in\{1,\dots,J\}$, belonging to the Dirichlet-Laplacian eigenbasis, for which, under the discretisation \eqref{Eq:Discretisation}, $\langle f,\phi_j\rangle_2=f_j$.
Accordingly, for $\bar{\mathbf f}_n$ and $\mathbf \Lambda_n$ as in \eqref{Eq:PostParam}, the plug-in posterior estimators are given by $\langle \bar f_n,\phi_j\rangle_2=\bar{\mathbf f}_{n,j}$, with corresponding posterior variances $\mathbf\Gamma_{n,jj}$. Thus, the $95\%$-credible interval for $\langle f,\phi_j\rangle_2$ is given by  $\bar{\mathbf f}_{n,j}\pm 1.96\sqrt{\mathbf\Gamma_{n,jj}}$.

	 In order to compute the asymptotic variances $\|\nabla\cdot(c\nabla \psi)\|_2^2$ appearing in the right hand side of \eqref{Eq:BvM} and \eqref{Eq:CLT}, we obtain the singular value decomposition (SVD) of the forward operator $G$, corresponding to finding the eigenfunctions $(\xi_k, \ k\in\N)\subset L^2(\Ocal)$ and eigenvalues $(\eta_k, \ k\in\N)\subset[0,\infty)$ solving the problem
\begin{equation}
\label{Eq:SVD}
\begin{split}
-\nabla\cdot(c\nabla \xi) - \eta \xi &=0, \ \ \textnormal{on}\ \ \Ocal \\
\xi&=0, \ \ \textnormal{on}\ \ \partial\Ocal,
\end{split}
\end{equation}
whereupon there follow the identities
$$
	G(f)=\sum_{k=1}^\infty \eta_k^{-1}\langle f, \xi_k\rangle_2\xi_k,
$$
and
$$
	\nabla\cdot(c\nabla u)=\sum_{k=1}^\infty \eta_k\langle u, \xi_k\rangle_2\xi_k;
	\qquad \|\nabla\cdot(c\nabla u)\|_2^2
	=\sum_{k=1}^\infty \eta_k^2\langle u, \xi_k\rangle_2^2.
$$

\begin{figure}[t!h]
\centering
\includegraphics[width=5.5cm,height=4.5cm]{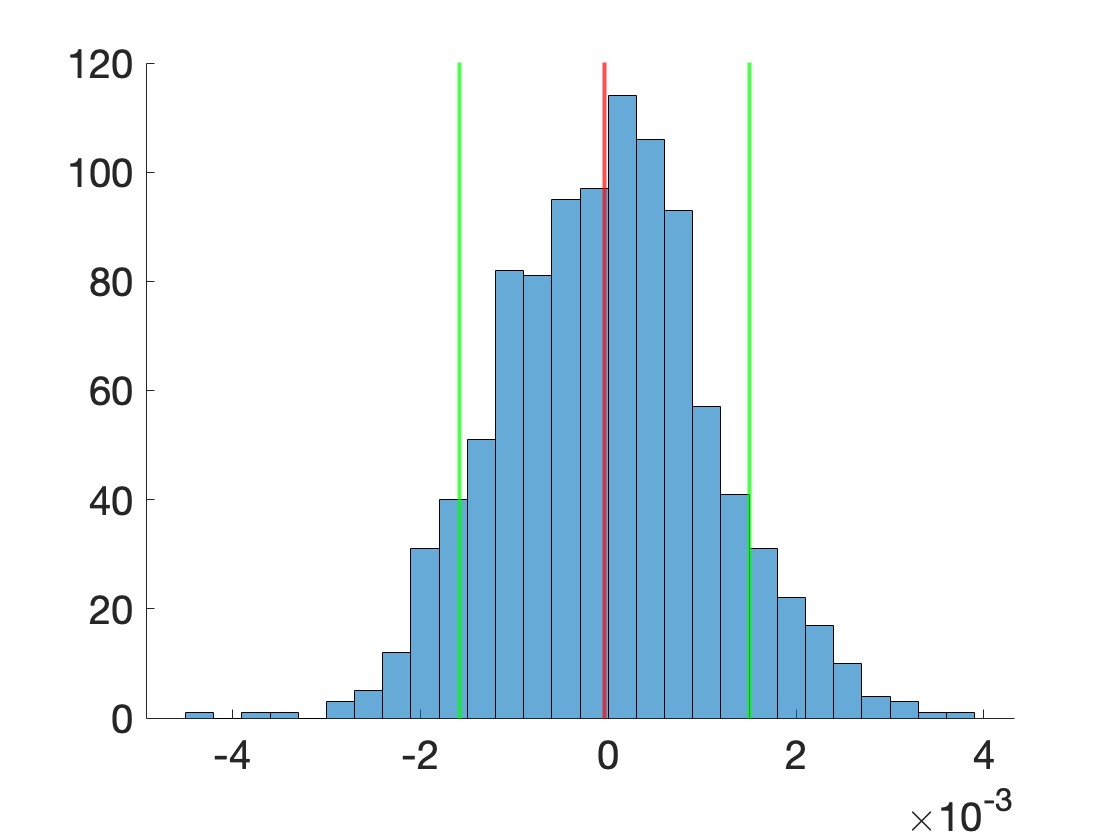}
\includegraphics[width=5.5cm,height=4.5cm]{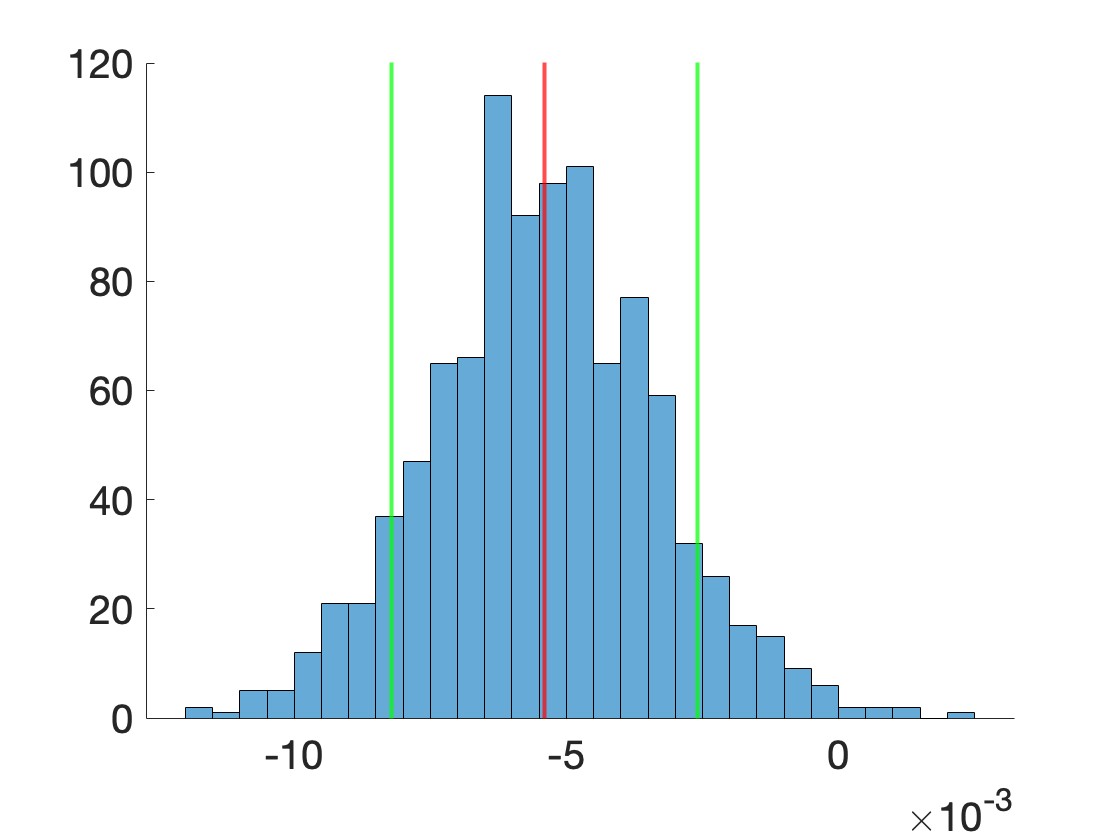}
\includegraphics[width=5.5cm,height=4.5cm]{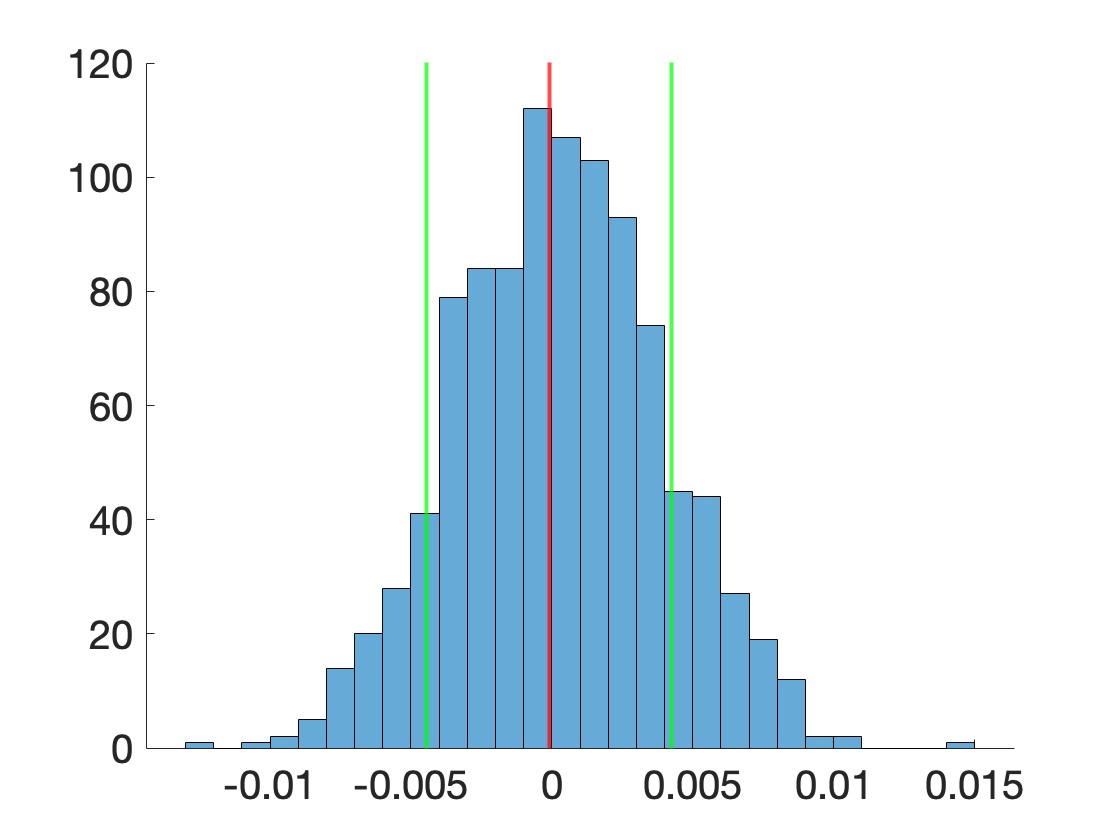}
\includegraphics[width=5.5cm,height=4.5cm]{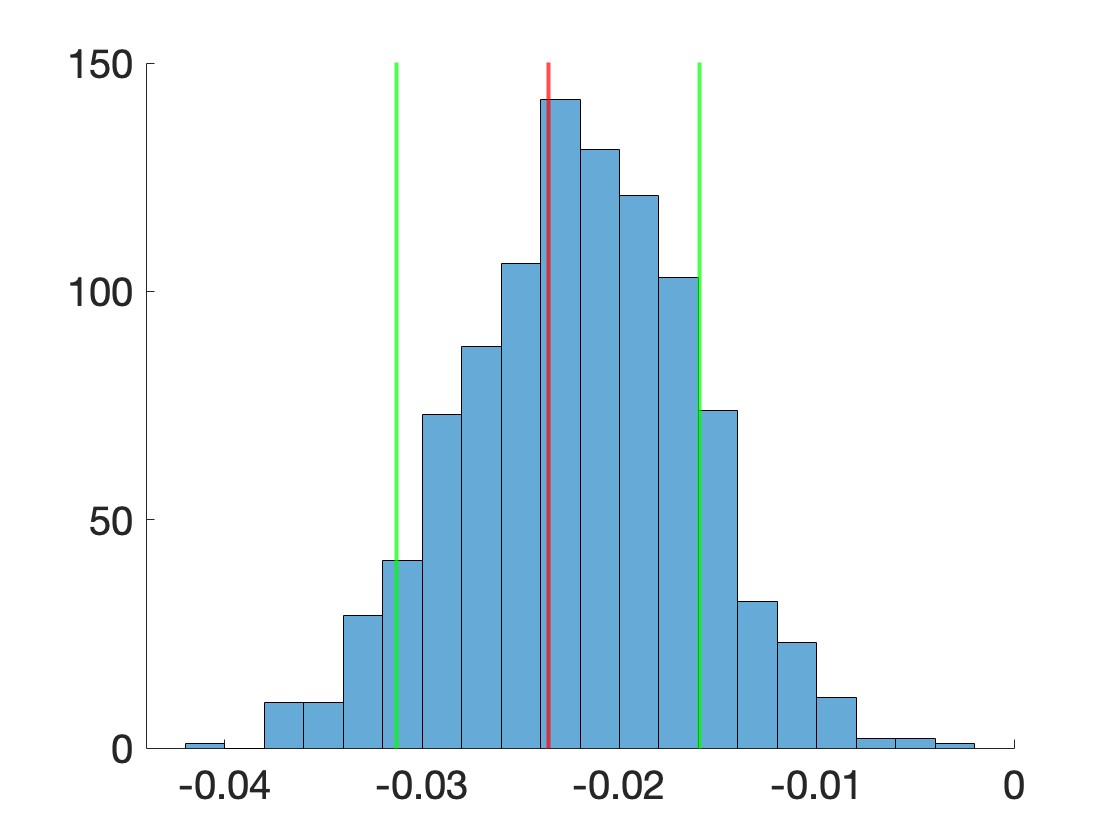}
\caption{Left to right, top to bottom: histograms relative to $1000$ realisations of the plug-in posterior mean estimators $\langle \bar f_n,\phi_j\rangle_2$, for $j=2,4,8,16$. The vertical red lines identify the true parameters $\langle f_0,\phi_j\rangle_2$. The vertical green lines identify the predicted asymptotic quantiles $ \langle f_0,\phi_j\rangle_2 \pm 1.96\|\nabla\cdot(c\nabla\phi_j)\|^2\sigma/\sqrt{n}$ (recalling the asymptotic regime $\varepsilon^2\simeq \sigma^2/n$). 
}
\label{Fig:Distributions}
\end{figure}

\begin{table}[!h]
\caption{Observed coverages for increasing sample sizes of the $95\%$-credible intervals for the linear functionals $\langle f,\phi_j\rangle_2$, with $j=2,4,8,16$.}
\label{Tab:Coverage} 
\begin{tabular}{p{2cm}p{1.75cm}p{1.75cm}p{1.75cm}p{1.75cm}p{1.75cm}p{1.75cm}}
\hline\noalign{\smallskip}
 $n$ & 50 & 100 & 250 & 500 & 750 & 1000   \\
\noalign{\smallskip}\hline\noalign{\smallskip}
$\phi_2$ & 0.885  & 0.921& 0.954 & 0.969 & 0.947 & 0.962 \\
\noalign{\smallskip}\hline\noalign{\smallskip}
$\phi_4$ & 0.904  & 0.936& 0.943 & 0.96 & 0.953 & 0.958  \\
\noalign{\smallskip}\hline\noalign{\smallskip}
$\phi_8$ & 0.92  & 0.934 & 0.951 &  0.944 & 0.963 & 0.952  \\
\noalign{\smallskip}\hline\noalign{\smallskip}
$\phi_{16}$ &  0.949 & 0.922 & 0.92 & 0.946 & 0.945  & 0.96  \\
\noalign{\smallskip}\hline\noalign{\smallskip}

\end{tabular}
\end{table}

In practice, we tackle the eigenvalue problem \eqref{Eq:SVD} via finite element methods exactly as outlined in Section \ref{Sec:Methodology} for the computation of the Dirichlet-Laplacian eigenbasis, obtaining numerical approximations of the eigenpairs. We note that, while used here as a convenient computational device to evaluate the asymptotic variances, knowledge of the SVD of the forward operator $G$ is not assumed for the theoretical results of Section \ref{Sec:Theory}, nor is required for the specification of the two classes of Gaussian priors introduced in Examples \ref{Ex:DirichLapl} and \ref{Ex:Matern} respectively. The theory and methodology investigated in the present article are indeed generally applicable to inverse problems where the SVD might be challenging or unfeasible to compute, or to settings where the properties of the associated eigenpairs might be unknown.

	Figure \ref{Fig:Distributions} shows the (approximate) distributions of the plug-in posterior mean estimators $\langle \bar f_n,\psi\rangle_2$ for four representative test functions $\psi=\phi_2,\phi_4,\phi_{8},\phi_{16}$. The plots present the histograms relative to $1000$ realisations of the estimators, obtained by drawing $1000$ independent collections of observations from the inverse regression model \eqref{Eq:DiscrModel}. For each experiment, a sample of size $n=1000$ was drawn, with noise standard deviation $\sigma=0.0005$. As expected from the central limit theorem \eqref{Eq:CLT}, the distributions of the plug-in estimators $\langle \bar f_n,\psi\rangle_2$ exhibit a normal shape, are approximately centred around the true parameter $\langle f_0,\psi\rangle_2$, and their spread is mostly captured by the asymptotic variance $\|\nabla\cdot(c\nabla \psi)\|^2_2$.

	Finally, Table \ref{Tab:Coverage} reports the coverage, for increasing sample sizes, of the $95\%$-credible intervals \eqref{Eq:CredInt} for the same linear functionals $\langle f,\phi_j\rangle_2$, with $j=2,4,8,16$, considered in the previous set of experiments. The results were obtained by drawing $1000$ independent collections, of size $n=1000$, of observations from the inverse regression model \eqref{Eq:DiscrModel}, with noise standard deviation $\sigma=0.0005$. For each random sample, a realisation of the $95\%$-credible intervals $\bar{\mathbf f}_{n,j}\pm 1.96\sqrt{\mathbf\Gamma_{n,jj}}$ was obtained, and the final coverage scores were computed as the fraction of times in which the true parameters $\langle f_0,\phi_j\rangle_2$ were contained in the obtained credible intervals.
As expected from the theoretical convergence result in \eqref{Eq:Coverage}, the observed coverages stabilise, as the sample size increases, around the prescribed credibility level $95\%$.

%
%
%

\subsection{Posterior inference with the Matérn process prior}
\label{Sec:Matern}

\begin{figure}[t]
\centering
\includegraphics[width=3.85cm,height=4cm]{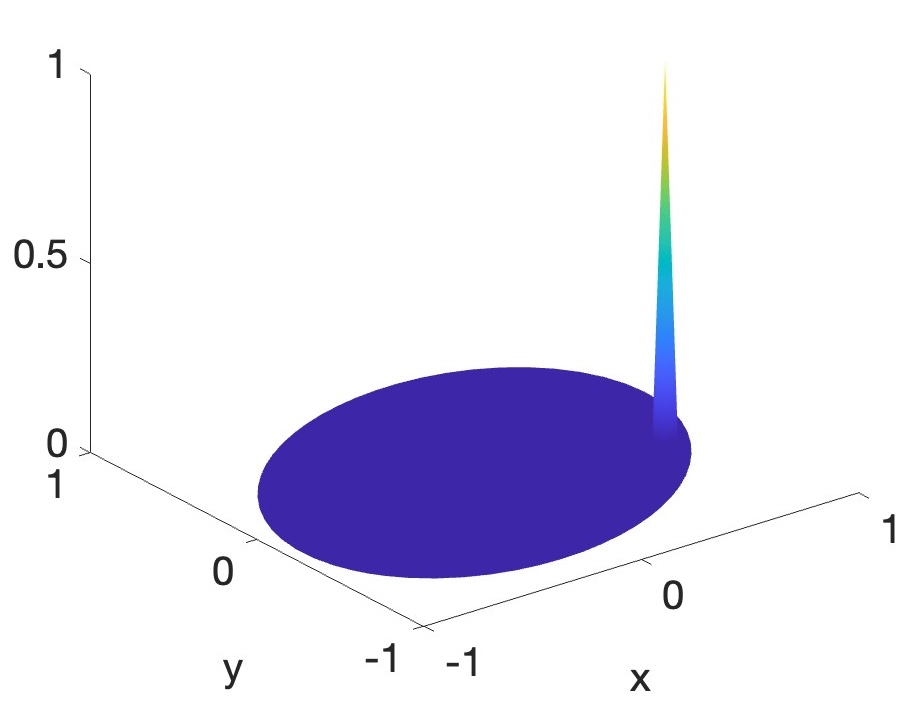}
\includegraphics[width=3.85cm,height=4cm]{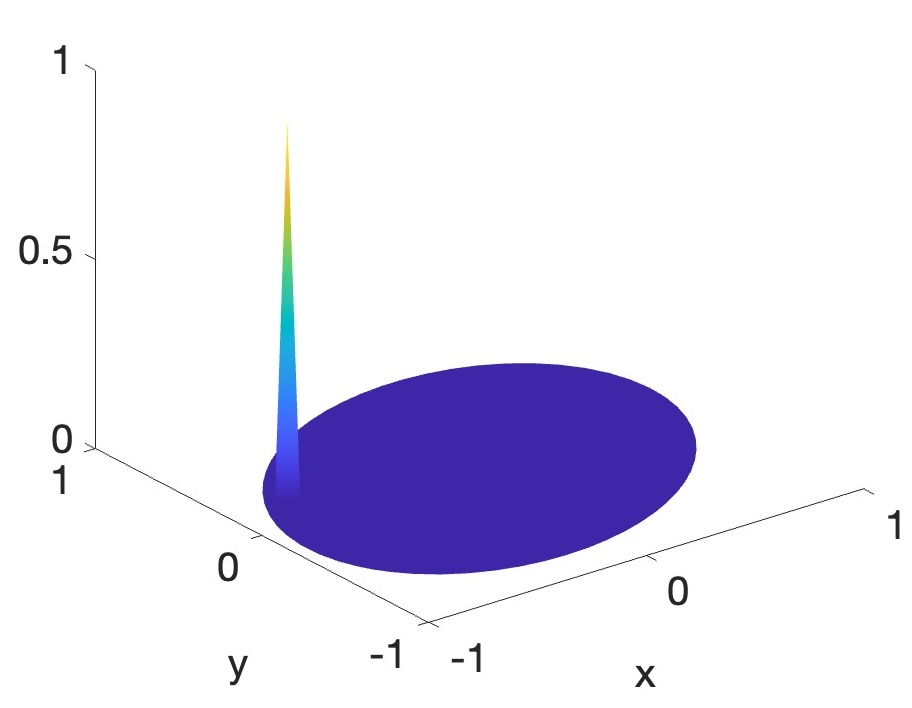}
\includegraphics[width=3.85cm,height=4cm]{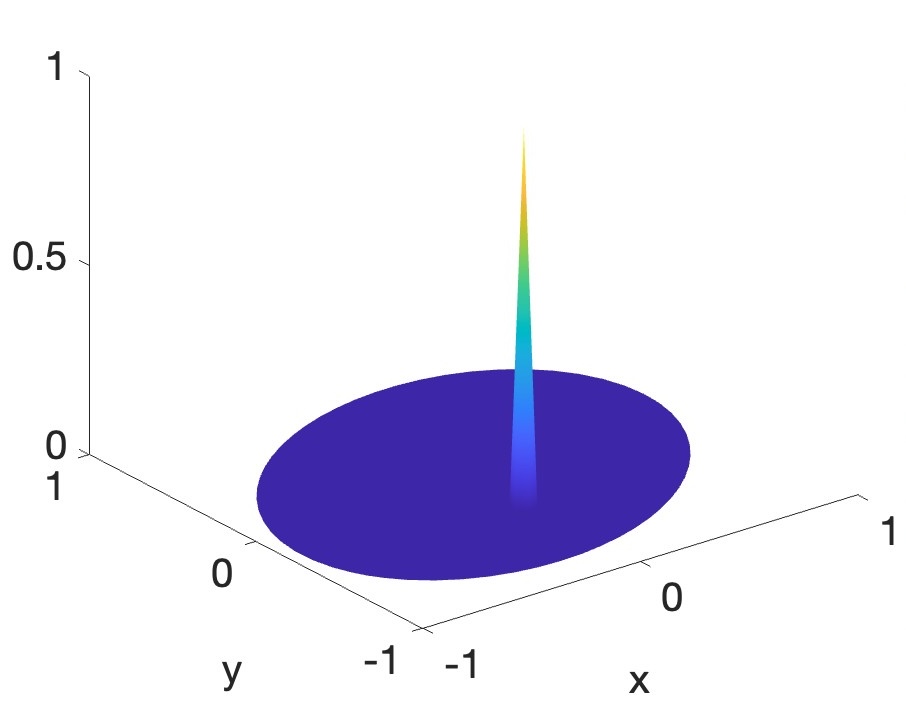}
\caption{Left to right: piecewise linear interpolation functions $\varphi_{75}$, $\varphi_{100}$ and $\varphi_{250}$ for a triangular mesh with $M=1169$ nodes. 
}
\label{Fig:LinInt}
\end{figure}

Next, we consider the Matérn process priors introduced in Example \ref{Ex:Matern}. We discretise the parameter space by assuming that $f$ is given by the finite sum
\begin{equation}
\label{Eq:Sum}
	f = \sum_{m=1}^M f_{m}\varphi_m, \qquad f_1,\dots,f_M\in\R,
	\qquad M\in\N,
\end{equation}
where $\varphi_1 ,\dots, \varphi_M$ are piecewise linear functions on the nodes $z_1,\dots,z_M\in\Ocal$ of a deterministic triangular mesh, uniquely identified by the property $\varphi_m(z_{m'}) = 1_{\{m=m'\}}$; see Figure \ref{Fig:LinInt}. Accordingly, $f$ in \eqref{Eq:Sum} satisfies $f(z_m) = f_m$, and for any $x\in\Ocal$ the value $f(x)$ is obtained by linearly interpolating the pairs $\{(z_m,f_m), \ m=1,\dots,M\}$.

	Under the discretisation \eqref{Eq:Sum}, the inverse regression model \eqref{Eq:DiscrModel} can be written in matrix notation exactly as in \eqref{Eq:MatrixModel}, now with
$$
	\mathbf G:=[G(\varphi_m)(x_i), \ i=1,\dots,n, \ m=1,\dots,M]\in\R^{n,M},
$$
and with
$$
	 \mathbf f:=(f_1,\dots,f_M)^T\in\R^M.
$$
Thus, again, $\mathbf Y|\mathbf f\sim N_n(\mathbf G\mathbf f,\sigma^2 \mathbf I_n)$. Similarly to Section \ref{Sec:Methodology}, the numerical computation of the matrix $\mathbf G$ can be carried out with finite element methods for elliptic PDEs.

	Recalling that, under the discretisation \eqref{Eq:Sum}, $f_m=f(z_m)$ for $m=1,\dots,M$, and the finite dimensional distributions property \eqref{Eq:FDDs}, assigning to $f$ a Matérn process prior with covariance $C_{\alpha,\ell}$ as in \eqref{Eq:MatCov} corresponds to assigning to the vector $\mathbf f$ the $M$-dimensional Gaussian prior with covariance matrix
$$
	\mathbf f \sim N_M(0,\mathbf C),
	\qquad \mathbf C:=[C_{\alpha,\ell}(z_h,z_m)]_{h,m=1}^M\in \R^{M,M}.
$$
The same conjugate computation as the one outlined in Section \ref{Sec:Methodology} can then be carried out, leading to the Gaussian posterior distribution $\mathbf f|\mathbf Y\sim N_M(\bar{\mathbf f}_n,\mathbf C_n)$, with posterior mean and covariance matrix respectively given by
$$
	\bar{\mathbf f}_n
	:=\frac{1}{\sigma^2}\mathbf C_n\mathbf G^T\mathbf Y; 
	\qquad \mathbf C_n:=(\sigma^{-2}\mathbf G^T\mathbf G+\mathbf C^{-1})^{-1}.
$$
Using the above conjugate formulae, posterior inference for the source function $f$ based on the Matérn process prior can efficiently be implemented.

\begin{figure}[t]
\centering
\includegraphics[width=5.5cm,height=4.5cm]{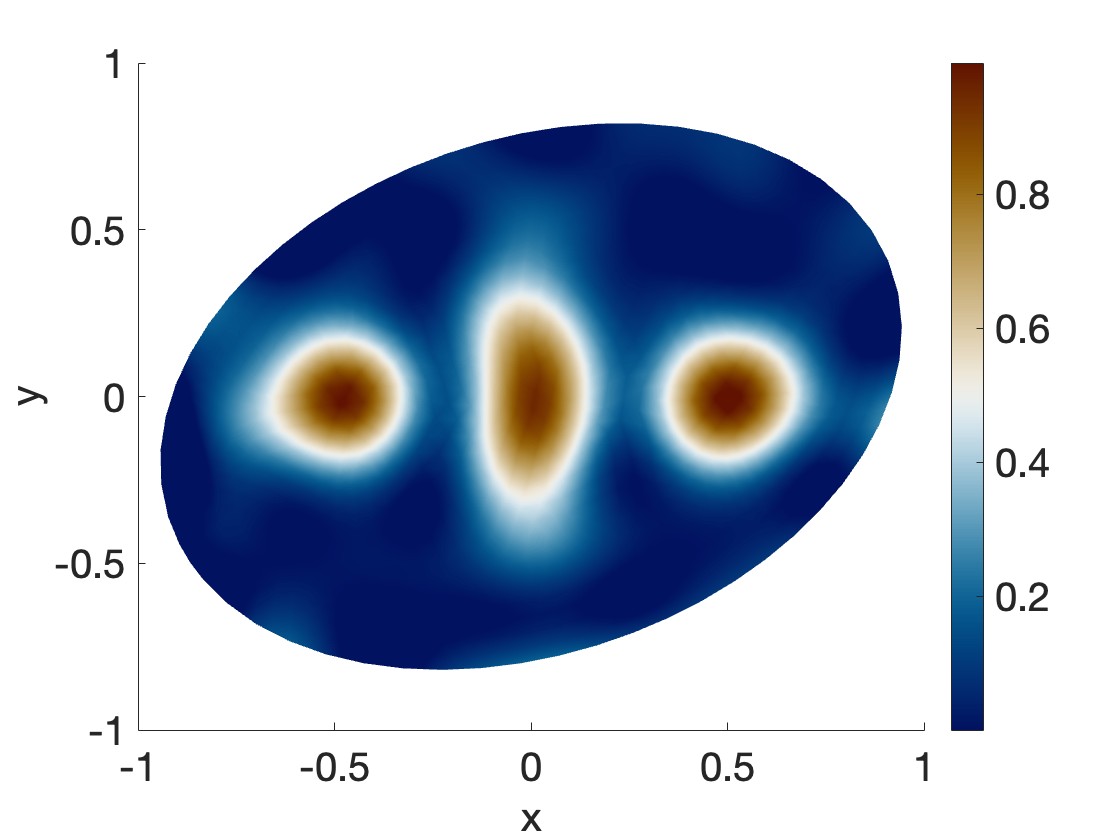}
\caption{The posterior mean estimate of the source function $f$ based on a Matérn process prior with parameters $\alpha=10$, $\ell=0.25$. 
}
\label{Fig:Matern}
\end{figure}

\begin{table}[t]
\caption{$L^2$-estimation errors achieved by the posterior mean estimator $\bar f_n$ arising from the Matérn process prior for increasing sample sizes.}
\label{Tab:MaternErrors} 
\begin{tabular}{p{3cm}p{.925cm}p{.925cm}p{.925cm}p{.925cm}p{.925cm}p{.925cm}p{.925cm}p{.925cm}p{.7cm}}
\hline\noalign{\smallskip}
 $n$ & 50 & 100 & 250 & 500 & 750 & 1000 &2000 & 3000 & 4500  \\
\noalign{\smallskip}\hline\noalign{\smallskip}
$\|\bar f_n - f_0\|_2$ & 0.30  & 0.30& 0.18 & 0.12 & 0.13& 0.10 & 0.086 & 0.076 & 0.067\\
\noalign{\smallskip}\hline\noalign{\smallskip}
$\|\bar f_n - f_0\|_2/\|f_0\|_2$ & 62.5\%  & 62.5\% & 37.5\%  & 25\% & 27.1\%&  20.8\% &17.9\%& 15.6\%& 13.9\% \\
\noalign{\smallskip}\hline\noalign{\smallskip}
\end{tabular}
\end{table}


%

	For the ground truth displayed in Figure \ref{Fig:StatProbl} (left), Table \ref{Tab:MaternErrors} reports the $L^2$-estimation errors attained by the posterior mean estimate $\bar f_n = \sum_{m=1}^M\bar{\mathbf f}_{n,m}\varphi_m$ based on an increasing number of observations from the inverse regression model \eqref{Eq:DiscrModel}, with noise standard deviation $\sigma=0.0005$ (with corresponding signal-to-noise ratio $\|G(f_0)\|_2/\sigma=37.55$). Across the experiments, the parameter space was discretised using a triangular mesh with $M=1169$ nodes. The prior regularity parameter for the Matérn process prior was set to $\alpha=10$, and the length-scale parameter to $\ell=0.25$. Figure \ref{Fig:Matern} shows the posterior mean estimate resulting from $n=4500$ observations.

	The results are relative to the same collection of synthetic data sets with increasing sample size employed in Section \ref{Sec:SeriesExp}, allowing a direct comparison to the results obtained with the Gaussian series priors considered therein. Overall, the achieved $L^2$-estimation errors are comparable in magnitude for each sample size, albeit the performance of the Gaussian series priors was consistently slightly better. It is plausible that such small discrepancy is caused by finite sample effects, prior tuning and the various numerical approximations.
	
%
%
%

\subsection{Further numerical experiments}

%
%
%

\subsubsection{Sensitivity to the noise variance}
\label{Sec:NoiseVar}

\begin{table}[t]
\caption{$L^2$-estimation errors achieved by the posterior mean estimator $\bar f_n$ arising from the Gaussian series prior for decreasing noise standard deviation.}
\label{Tab:sigmaErrors} 
\begin{tabular}{p{3cm}p{1.75cm}p{1.75cm}p{1.75cm}p{1.75cm}p{1.5cm}p{.9cm}}
\hline\noalign{\smallskip}
 $\sigma$ & 0.01  & 0.005  & 0.0025 & 0.001 &  0.0005& 0.0001  \\
 \noalign{\smallskip}\noalign{\smallskip}
$\|G(f_0)\|_2/\sigma$ & 1.88  &3.75 & 7.51 & 18.77  & 37.55& 187.74 \\
\noalign{\smallskip}\hline\noalign{\smallskip}
$\|\bar f_n - f_0\|_2$ & 0.21  & 0.16 & 0.14  & 0.078 & 0.060&  0.049 \\
\noalign{\smallskip}\hline\noalign{\smallskip}
$\|\bar f_n - f_0\|_2/\|f_0\|_2$ & 43.75\%  & 33.33\% & 29.17\%  & 16.25\% & 12.5\%&  10.20\% \\
\noalign{\smallskip}\hline\noalign{\smallskip}
\end{tabular}
\end{table}

For the empirical results presented in Sections \ref{Sec:SeriesPriors} and \ref{Sec:Matern}, a fixed noise standard deviation $\sigma = 0.0005$ in the inverse regression model \eqref{Eq:DiscrModel} was used (with corresponding signal-to-noise ratio $\|G(f_0)\|_2/\sigma=37.55$). Here, we provide a brief investigation of the sensitivity of the considered methodology to the value of $\sigma$, performing a set of experiments with decreasing standard deviation from $\sigma=0.01$ (for which $\|G(f_0)\|_2/\sigma=1.8774$) to $\sigma=0.0001$ (for which$\|G(f_0)\|_2/\sigma=187.74$), based on the same domain and ground truth used previously. Across the experiments, the sample size was kept fixed at $n=4500$.

For concreteness, we focus on the Gaussian series priors from Section \ref{Sec:SeriesPriors}, with the same prior tuning employed therein; similar results may be obtained with the Matérn process priors. The $L^2$-estimation error associated to the resulting posterior mean estimates are shown in Table \ref{Tab:sigmaErrors}. Unsurprisingly, these were observed to decrease monotonically as the signal-to-noise ratio increased.  In particular, at the lowest value $\sigma=0.0001$, the $L^2$-estimation error may be seen to approach the  $L^2$-approximation error resulting from projecting $f_0$ onto the employed basis, which is equal to $0.0486$.

%

\subsubsection{Inference with unknown noise variance}
\label{Sec:VarEstim}

We conclude the simulation study considering the important practical scenario where the noise standard deviation $\sigma$ is itself unknown and needs to be estimated from the data. Given observations $\{(Y_i,X_i)\}_{i=1}^n$ from the inverse regression model \eqref{Eq:DiscrModel}, we undertake the simple `empirical Bayes' approach of obtaining a preliminary estimate $\hat \sigma_n$ of $\sigma$, and then carry over the methodology laid out in Sections \ref{Sec:SeriesPriors} and \ref{Sec:Matern} with $\sigma$ replaced by $\hat\sigma_n$. Alternatively, we note that a joint Bayesian model for $f$ and $\sigma$ in \eqref{Eq:DiscrModel} could be considered by endowing $\sigma$ with a prior distribution. For example, an independent inverse-gamma distribution would lead (conditionally given $f$) to a conjugate posterior distribution, whereupon joint posterior samples for the pair $(f,\sigma)$ could readily be obtained via a Gibbs sampler, alternating draws from the Gaussian posterior distribution of $f|\{(Y_i,X_i)\}_{i=1}^n,\sigma$ and draws from the inverse-gamma posterior distribution of $\sigma|\{(Y_i,X_i)\}_{I=1}^n,f$. For brevity, we will not pursue this approach further here.

\begin{table}[t]
\caption{Inferential results for the difference-based estimator $\hat\sigma_n$ of the noise standard deviation and the the empirical Bayes posterior mean estimator $\hat f_n$, for increasing sample sizes.}
\label{Tab:sigma2} 
\begin{tabular}{p{2.75cm}p{2.75cm}p{2.75cm}p{2.75cm}p{1.92cm}}
\hline\noalign{\smallskip}
$n$ & 1000  &2000 & 3000 & 4500  \\
\noalign{\smallskip}\hline\noalign{\smallskip}
$\hat \sigma $ & 0.0034  & 0.0024 & 0.0023  & 0.00072 \\
\noalign{\smallskip}\hline\noalign{\smallskip}
$\|\hat f_n - f_0\|_2$ & 0.18  & 0.14 & 0.12  & 0.063  \\
\noalign{\smallskip}\hline\noalign{\smallskip}
$\|\hat f_n - f_0\|_2/\|f_0\|_2$ & 37.5\%  & 29.17\% & 25\%  & 13.12\%  \\
\noalign{\smallskip}\hline\noalign{\smallskip}
\end{tabular}
\end{table}

	Several strategies have been proposed in the literature for variance estimation in nonparametric regression models, ranging from residual-based estimators using kernel smoothing \cite{hall1990variance} and splines \cite{wahba1978improper}, to difference-based estimators \cite{rice1984bandwidth}. See \cite{alharbi2013error} for an overview. Here, we will consider the difference-based method proposed in \cite{rice1984bandwidth}, estimating $\sigma$ in model \eqref{Eq:DiscrModel} by
$$
	\hat\sigma_n := \sqrt{\hat\sigma^2_n},
	\qquad \hat\sigma^2_n:=\frac{1}{2(n-1)}\sum_{i=2}^n (Y_i - Y_{i-1})^2.
$$
Based on $\hat\sigma_n$, the `empirical Bayes posterior mean' estimate $\hat f_n$ arising from a Gaussian series prior or a Matérn process prior can then be readily computed exactly as described in Sections \ref{Sec:SeriesPriors} and \ref{Sec:Matern} respectively, replacing $\sigma$ with $\hat \sigma_n$ in the relevant conjugate formulae.

Table \ref{Tab:sigma2} summarises the inferential results obtained with the  difference-based estimation procedure for increasing sample sizes. For these experiments, the noise standard deviation was set to $\sigma=0.0005$, and a Gaussian series priors with the same tuning as in Sections \ref{Sec:SeriesExp} and \ref{Sec:NoiseVar} was used. The results show a progressive improvement in the reconstruction quality for both the noise standard deviation $\sigma$ and the unknown source function $f$. In particular, for the largest considered sample size $n=4500$, the $L^2$-estimation error $\|\hat f_n - f_0\|_2$ resulted to be only marginally higher than the one obtained under the same experimental conditions (and prior tuning) in the context of the empirical results presented in Section \ref{Sec:SeriesExp} (cf.~Table \ref{Tab:Errors}), for which knowledge of the value of $\sigma$ was assumed.

%

%
%
%
%
%

\section{Summary and discussion}

In this article we have considered the nonparametric Bayesian approach with Gaussian priors to linear inverse problems, focusing on the important example of source identification in elliptic PDEs. The main advantages of the considered methodology lie in its modelling flexibility, its ease of implementation (cf.~the conjugate formulae \eqref{Eq:Conjugacy} and \eqref{Eq:PostParam}), as well as its theoretical guarantees on estimation and uncertainty quantification (cf.~Section \ref{Sec:Theory}). The performance of the approach has been investigated in a numerical simulation study (cf.~Section \ref{Sec:Numerics}) under two distinct prior models (Gaussian series and Matérn process priors), both for which excellent reconstruction results have been obtained.

	The present work also raises various related research questions. Firstly, it is of interest and practical importance to further explore the setting where the noise standard deviation $\sigma$ in the inverse regression model \eqref{Eq:DiscrModel} is unknown. While the simple difference-based estimator considered in Section \ref{Sec:VarEstim} has proved effective, several competing approaches, including the joint conjugate Bayesian model outlined in Section \ref{Sec:VarEstim}, could be investigated. Furthermore, a related interesting question concerns the extensions of the theoretical results presented in Section \ref{Sec:Theory} to the setting with unknown variance; see e.g.~\cite{kejzlar2021fast} for related results in a direct regression model.

	Lastly, let us mention the important issue of specifying the hyperparameter values for the considered prior distributions,  namely the truncation level and regularity in the Gaussian series \eqref{Eq:DiscrPrior}, and the smoothness and length-scale parameters in the Matérn covariance kernel \eqref{Eq:MatCov}. There is by now a vast literature investigating the methodological and theoretical aspects of empirical and hierarchical Bayesian approaches to fully data-driven selection of the hyperparameters; see \cite{KSvdVvZ15, roussszabo,teckentrup2020convergence,agapiou2014analysis} and the many reference therein. Investigating the implications and performance of these methods in the context of the observation model and prior distributions considered in the present article is an interesting problem for future research. 

%
%
%
%
%

\paragraph{Acknowledgement.} The Author is grateful to three anonymous referee for many helpful comments that lead to an improvement of the article. This research has been partially supported by MUR, PRIN project
2022CLTYP4. The Author also gratefully acknowledges the support of ``de Castro" Statistics Initiative, Collegio Carlo Alberto, Torino. There are no conflicts of interest to declare that are relevant to the content of this chapter.

\bibliographystyle{acm}
\bibliography{LinIPRef}

\end{document}